\documentclass[12pt]{article}
\usepackage{amsmath}    
\usepackage{amssymb}    
\usepackage{amsfonts}   
\usepackage{amsthm}     
\usepackage{geometry}   
\usepackage{graphicx}   
\usepackage{enumitem}   
\usepackage{indentfirst}
\usepackage{tikz-cd}    
\usepackage{caption}    
\usepackage{mathrsfs}   
\usepackage[normalem]{ulem}  
\usepackage{color}      
\usepackage{mathtools}  
\usepackage{tcolorbox}  
\usepackage[hidelinks]{hyperref}
\usepackage{doi}

\usepackage[all]{xy}

\usepackage[numbers]{natbib}
\providecommand{\keywords}[1]{\textbf{Keywords:} #1}
\providecommand{\subjclass}[1]{\textbf{2020 Mathematics Subject Classification:} #1}
\newtheorem{definition}{Definition}[section]
\newtheorem{theorem}[definition]{Theorem}
\newtheorem{proposition}[definition]{Proposition}
\newtheorem{lemma}[definition]{Lemma}
\newtheorem{example}[definition]{Example}
\newtheorem{remark}[definition]{Remark}

\newtheorem{corollary}[definition]{Corollary}
\numberwithin{equation}{section}

\begin{document}
	
	\title{Topological Mundici Adjunction and Its Applications
}
	\author{Yanqiao Jia
		, Xueting Zhao
		, Yichuan Yang\thanks{Corresponding author. E-mail: ycyang@buaa.edu.cn}}
	\date{}
	
	\maketitle

\begin{abstract}
	We establish a canonical topological lifting of Mundici equivalence, representing every topological MV-algebra as the unit interval of a topological $\ell$-group and extending the correspondence to a categorical adjunction.
	As applications, we give a direct group-theoretic proof of the equivalence of $T_0$ and $T_{3\frac12}$ and show that every topological MV-algebra is a retral space.
	We further extend the lifting to topological Riesz MV-algebras and identify the canonical topology with Wiweger's mixed topology under suitable assumptions, recovering the strict topology as a special case.
	Finally, we obtain a unit-interval formulation of Koopman representations for continuous-time dynamical systems.
\end{abstract}
	
	
	\medskip

\section*{Introduction}\label{sec:introduction}

MV-algebras were introduced by Chang~\cite{chang1958} as an algebraic framework for infinite-valued {\L}ukasiewicz logic. A fundamental development in their structural theory was Mundici equivalence~\cite{mundici1986}, which establishes a categorical equivalence between MV-algebras and abelian lattice-ordered groups ($\ell$-groups for short) with a distinguished strong unit. 
This equivalence has played a central role in the development
of MV-algebra theory and has inspired extensive research
across several areas of mathematics
(see, e.g., \cite{MarraReggio2017,MarraMenni2024,Mundici2018,Wehrung2019}).

A further development is the theory of Riesz MV-algebras, which extend MV-algebras with scalar multiplication by elements of $[0,1]$ and are categorically equivalent to Riesz spaces with strong unit~\cite{DiNolaFlondorLeustean2003,DiNolaLeustean2014}. This development has further strengthened the connection between MV-algebra theory and ordered functional analysis (see, e.g., \cite{BedeDiNola2004,DiNolaLeustean2011,DiNolaLapentaLeustean2018}).

On the topological side, both topological $\ell$-groups
(see, e.g., \cite{Ball1979,Redfield1974,Jordan2023})
and topological Riesz spaces
(see, e.g., \cite{Aliprantis1974,AliprantisBurkinshaw2003})
have a long history of research.
The topological study of MV-algebras dates back to Hoo's pioneering work~\cite{hoo1997}. Subsequent research has revealed that topological MV-algebras, although their additive reducts are topological monoids rather than groups, share remarkable topological properties with topological groups (see, e.g., \cite{gan2024,li2026,Weber2012}).

This resemblance suggests a natural connection between
topological MV-algebras and topological $\ell$-groups.
Recently, Li, Mu, and Yang~\cite{LiMuYang2027}
explicitly raised the question of how topological
MV-algebras can be embedded into topological $\ell$-groups,
or even into topological Riesz spaces.
In view of Mundici equivalence, we formulate the following
more precise \emph{topological $\ell$-group lifting question}:

\medskip
\emph{Is every topological MV-algebra
	topologically isomorphic to the unit interval
	of a topological $\ell$-group?}
\medskip

The analogous Riesz-space lifting question is:

\medskip
\emph{Is every topological Riesz MV-algebra
	topologically isomorphic to the unit interval
	of a topological Riesz space?}
\medskip

Building on Ball's earlier work on convergence and Cauchy structures in $\ell$-groups~\cite{Ball1980}, Ball, Georgescu, and Leu\c{s}tean~\cite{BallGeorgescuLeustean2002} established a correspondence between MV-convergences and bounded $\ell$-convergences on the associated Mundici group.
Since topologies can be characterized by convergence, their approach suggested a possible solution to the topological lifting problem. 
However, the assertion concerning the preservation of
topologicity in~\cite[Remark~3.4(2)]{BallGeorgescuLeustean2002}
does not hold in general, as demonstrated by the counterexample
in Remark~\ref{rem:canonical-topology}(2).
Thus, this correspondence alone does not establish
the desired topological lifting.

In earlier work~\cite{JiaYang2026}, we investigated this problem through the notion of \emph{summability}. We proved that a topological MV-algebra can be represented as the unit interval of a lattice-ordered topological group if and only if it is summable, and introduced \emph{strong summability} as the corresponding criterion for representations in topological $\ell$-groups. We also established categorical adjunctions extending Mundici equivalence to the corresponding classes of topological structures. However, whether these conditions are automatically satisfied by arbitrary topological MV-algebras remained open.

\medskip
\noindent\textbf{The present paper.}
We resolve the topological $\ell$-group lifting question
raised above in full generality.
We prove that \emph{every topological MV-algebra is
	strongly summable}.
Consequently, the topology of every topological
MV-algebra admits a canonical lifting to its Mundici
group, making the latter a topological $\ell$-group
whose unit interval recovers the original topological
MV-algebra.
No separation, completeness, or local convexity
assumptions are required.
This removes the restrictions in the preceding
representation criteria and extends the categorical
topological lifting of Mundici equivalence to all
topological MV-algebras.

The significance of this construction goes beyond the existence of a compatible group topology. Since the lifted group is a topological $\ell$-group, its lattice operations provide a canonical continuous retraction onto the unit interval. This connects topological MV-algebras with the theory of retral spaces and allows classical results on topological groups and covering properties to be transferred to the MV-algebra setting. Moreover, the construction extends naturally to topological Riesz MV-algebras and their associated Riesz spaces. Under suitable hypotheses, we show that the resulting topology coincides with Wiweger's mixed topology~\cite{Wiweger1961}, thereby connecting the topological lifting problem with classical functional analysis. This connection further leads to applications involving Koopman representations of dynamical systems.

\medskip
\noindent\textbf{The main results of this paper are as follows.}

\begin{itemize}
	\item \emph{Topological lifting of Mundici equivalence.}
	We prove that every topological MV-algebra can be
	represented as the unit interval of a topological
	$\ell$-group (Theorem~\ref{thm:topological-lifting})
	and establish the corresponding topological Mundici
	adjunction (Theorem~\ref{thm:categorical-topological-lifting}).
	Together, these results provide a complete affirmative
	answer to the topological $\ell$-group lifting question
	formulated above, at both the representational and
	categorical levels.

	\item \emph{Topological and structural applications.}
	We give a short and direct group-theoretic proof of the equivalence of $T_0$ and $T_{3\frac12}$ (Corollary~\ref{cor:separation-equivalence}), a previously open separation problem recently resolved
	in~\cite{Luan2026,XieYangUniformization2026}. 
	We also prove that every topological MV-algebra is a retral space (Proposition~\ref{prop:canonical-retract}), thereby providing an affirmative answer to the question posed by Feng and Gartside~\cite[Question~4.16(i)]{FengGartside2026} within this class.
	
	\item \emph{Topological lifting of Riesz MV-algebras.}
	We extend the canonical topological lifting to Riesz
	MV-algebras and establish a corresponding adjunction
	between topological Riesz MV-algebras and topological
	Riesz spaces with strong unit
	(Proposition~\ref{prop:topological-riesz-lifting}),
	thereby also providing an affirmative answer to the
	Riesz-space lifting question formulated above.
	
	\item \emph{Functional-analytic and dynamical applications.}
	We show that, under suitable assumptions, the canonical
	lifted topology coincides with Wiweger's mixed topology,
	recovering the strict topology on $C_b(X)$ as a special
	case (Proposition~\ref{prop:canonical-mixedification}).
	We also obtain a unit-interval formulation of Koopman
	representations for continuous-time dynamical systems
	(Theorem~\ref{thm:unit-interval-koopman}).
\end{itemize}

Taken together, these results establish a topological Mundici adjunction and forge new connections between the theory of topological MV-algebras, functional analysis, and dynamical systems.

The paper is organized as follows. Section~\ref{sec:preliminaries} recalls the necessary preliminaries, and Section~\ref{sec:Topological_lifting} establishes the topological lifting of Mundici equivalence. Section~\ref{sec:Applications} presents applications: Subsection~\ref{subsec:topological-properties} studies topological properties of MV-algebras, while Subsection~\ref{subsec:riesz-koopman} develops the lifting theory for Riesz MV-algebras and its applications to functional analysis and dynamical systems. Finally, Section~\ref{sec:conclusion} concludes the paper and discusses future research directions.

\medskip
\noindent\textbf{Note added.}
Almost simultaneously, Meng, Shen, and
Zhang~\cite{MengShenZhang2026}
independently established a topological
Mundici representation theorem and the corresponding
adjunction using a different approach.

\section{Preliminaries}\label{sec:preliminaries}
In this section, we recall the necessary background
and fix notation used throughout the paper.

We begin by recalling some basic notions concerning
MV-algebras and lattice-ordered groups.

\begin{definition}[\cite{chang1958,mundici2000}]
	An \emph{MV-algebra} is an algebra $(A,\oplus,\neg,0)$
	of type $(2,1,0)$ satisfying the following axioms
	for all $x,y\in A$:
	\begin{enumerate}[label=(M\arabic*)]
		\item $(A,\oplus,0)$ is a commutative monoid;
		\item $\neg\neg x=x$;
		\item $x\oplus\neg0=\neg0$;
		\item $\neg(\neg x\oplus y)\oplus y
		=\neg(\neg y\oplus x)\oplus x$.
	\end{enumerate}
\end{definition}

Set $1:=\neg0$ and define the derived operations
\[
x\odot y:=\neg(\neg x\oplus\neg y),
\qquad
x\ominus y:=x\odot\neg y.
\]
The following elementary identities hold in every
MV-algebra:

\[
\begin{gathered}
	\neg1=0,\qquad x\oplus1=1,\qquad
	x\oplus\neg x=1,\\
	x\odot\neg x=0,\qquad
	x\oplus y=\neg(\neg x\odot\neg y).
\end{gathered}
\]

The \emph{natural order} on $A$ is defined by
$x\leq y$ if and only if $\neg x\oplus y=1$.
With this order, $A$ is a bounded distributive
lattice with bottom element $0$, top element $1$,
and lattice operations
\[
x\vee y=y\oplus(x\odot\neg y),
\qquad
x\wedge y=x\odot(\neg x\oplus y).
\]

An \emph{MV-homomorphism} is a map preserving
$\oplus$, $\neg$, and $0$. We denote by
\[
\mathcal{MV}
\]
the category of MV-algebras and MV-homomorphisms.

\begin{example}
	The unit interval $[0,1]$, equipped with
	\[
	x\oplus y=\min\{x+y,1\},
	\qquad
	\neg x=1-x,
	\]
	is an MV-algebra, called the
	\emph{standard MV-algebra}.
	Its derived operations are
	\[
	x\odot y=\max\{x+y-1,0\},
	\qquad
	x\ominus y=\max\{x-y,0\}.
	\]
\end{example}

We next recall a basic fact about MV-terms
that will be used later.
An \emph{MV-term} in variables $x_1,x_2,\ldots,x_n$
is an expression obtained from these variables
and the constant $0$ by finitely many applications
of $\oplus$ and $\neg$.
Every MV-term induces a function $A^n\to A$
on any MV-algebra $A$, called an
\emph{MV-term function};
see \cite[Section~1.4]{mundici2000}.

\begin{lemma}[{\cite[Lemma~3.1.9]{mundici2000}}]
	\label{lem:integer-affine-mv-term}
	Let $g:[0,1]^n\to\mathbb R$ be an
	integer-affine function, that is,
	\[
	g(x_1,x_2,\ldots,x_n)
	=b+\sum_{i=1}^n m_i x_i,
	\qquad b,m_1,m_2,\ldots,m_n\in\mathbb Z.
	\]
	Then $g^\sharp:=(g\vee0)\wedge1$
	is represented by an MV-term.
\end{lemma}

We now recall some basic notions and results
from the theory of lattice-ordered groups.
For a comprehensive treatment, see
\cite{birkhoff1967,darnel1994}.

A \emph{lattice-ordered group} (or \emph{$\ell$-group})
is a group $(G,+,0)$ equipped with a lattice order
that is invariant under both left and right translations.
An $\ell$-group is called \emph{abelian} if its underlying
group is commutative.
Throughout this paper, all $\ell$-groups are assumed
to be abelian. For $g\in G$, put
\[
g^+=g\vee0,\qquad
g^-=(-g)\vee0,\qquad
|g|=g^++g^-.
\]

An \emph{$\ell$-group homomorphism} is a group
homomorphism preserving finite joins and meets.
A positive element $u\in G$ is called a
\emph{strong unit} if for every $g\in G$
there exists $n\in\mathbb N$ such that
$|g|\leq nu$. We write $(G,u)$ for an
abelian $\ell$-group with a distinguished
strong unit.

We recall the following basic properties of
abelian $\ell$-groups.

\begin{proposition}[{\cite{birkhoff1967,darnel1994}}]
	\label{prop:subdirect-representation}
	Let $G$ be an abelian $\ell$-group and $x,y,z\in G$.
	The following hold:
	\begin{enumerate}[label=(\arabic*)]
		\item $G$ is a distributive lattice,
		$x+(y\vee z)=(x+y)\vee(x+z)$ and
		$x+(y\wedge z)=(x+y)\wedge(x+z)$;
		
		\item $-(x\vee y)=(-x)\wedge(-y)$ and
		$-(x\wedge y)=(-x)\vee(-y)$;
		
		\item $x=x^+-x^-$,
		$|x|=x\vee(-x)=x^+\vee x^-$,
		$x^+\wedge x^-=0$;
		
		\item $G$ is $\ell$-isomorphic to a
		subdirect product of totally ordered
		abelian $\ell$-groups. Consequently,
		an identity involving group and lattice
		operations holds in every abelian
		$\ell$-group if it holds in every
		totally ordered abelian $\ell$-group.
	\end{enumerate}
\end{proposition}

An abelian $\ell$-group gives rise to
an MV-algebra through the following construction.

\begin{example}[\cite{mundici2000}]
	\label{ex:lgroup-to-MV}
	Let $G$ be an abelian $\ell$-group and $u\in G$
	a positive element. Then
	$[0,u]=\{g\in G:0\leq g\leq u\}$
	is an MV-algebra with operations
	\[
	x\oplus y=(x+y)\wedge u,
	\qquad \neg x=u-x.
	\]
	In this case, $x\odot y=(x+y-u)\vee0$
	and $x\ominus y=(x-y)^+$.
\end{example}

Let
\[
\ell\mathcal G_u
\]
denote the category of abelian $\ell$-groups with
distinguished strong units and unit-preserving
$\ell$-group homomorphisms.
The construction in Example~\ref{ex:lgroup-to-MV}
defines a functor
\[
\Gamma:\ell\mathcal G_u\longrightarrow\mathcal{MV},
\qquad (G,u)\longmapsto[0,u],
\]
which sends each morphism $f:(G,u)\to(H,v)$
to $\Gamma(f)=f|_{[0,u]}$.

In fact, every MV-algebra is isomorphic to
one obtained by this construction, as established
by the following theorem.

\begin{theorem}[Mundici equivalence, \cite{mundici1986}]
	\label{thm:mundici}
	The functor $\Gamma:\ell\mathcal G_u\to\mathcal{MV}$
	is an equivalence of categories.
\end{theorem}

In particular, there exists a quasi-inverse functor
$\Xi:\mathcal{MV}\to\ell\mathcal G_u$ and a natural
isomorphism
$\varphi:\operatorname{id}_{\mathcal{MV}}\to\Gamma\Xi$.

For an MV-algebra $A$, write
$\Xi(A)=(G_A,u_A)$ and $\Gamma(G,u)=[0,u]$
when no confusion arises.
We call $(G_A,u_A)$ the \emph{Mundici group} of $A$.
The canonical MV-isomorphism
\[
\varphi_A:A\longrightarrow[0,u_A]
\]
satisfies
$\varphi_A(\neg x)=u_A-\varphi_A(x)$ and
$\varphi_A(x\oplus y)
=(\varphi_A(x)+\varphi_A(y))\wedge u_A$
for all $x,y\in A$.
For further details, see \cite[pp.~139--149]{mundici2000}.

We next recall the Riesz-space counterpart of
Mundici equivalence.

\begin{definition}
	\label{def:riesz-space}
	A \emph{Riesz space} (or \emph{vector lattice})
	is a real vector space $E$ whose underlying
	ordered additive group is an abelian $\ell$-group,
	such that $rx\geq0$ whenever $r\geq0$ and $x\geq0$.
\end{definition}

The notions of strong unit, positive part,
negative part, and modulus are inherited
from the underlying $\ell$-group.
A \emph{Riesz homomorphism} is a real-linear
map preserving finite joins and meets.

\begin{definition}[{\cite[Definition~3.1]{DiNolaLeustean2014}}]
	\label{def:riesz-mv-algebra}
	A \emph{Riesz MV-algebra} is an MV-algebra $A$
	equipped with a scalar multiplication
	$[0,1]\times A\to A$, $(r,x)\mapsto rx$,
	satisfying the following identities for all
	$r,q\in[0,1]$ and $x,y\in A$:
	\begin{align*}
		r(x\ominus y)&=rx\ominus ry,\\
		(r\ominus q)x&=rx\ominus qx,\\
		r(qx)&=(rq)x,\\
		1x&=x,
	\end{align*}
	where $r\ominus q=\max\{r-q,0\}$.
\end{definition}

A \emph{Riesz MV-homomorphism} is an
MV-homomorphism preserving scalar multiplication.
We denote by
\[
\mathcal{RMV}
\]
the category of Riesz MV-algebras and Riesz
MV-homomorphisms, and by
\[
u\mathcal R\mathit{iesz}
\]
the category of Riesz spaces with distinguished
strong unit and unit-preserving Riesz homomorphisms.

For $(E,u)\in u\mathcal R\mathit{iesz}$,
the interval $[0,u]$ is a Riesz MV-algebra,
with MV-operations as in
Example~\ref{ex:lgroup-to-MV} and scalar
multiplication inherited from $E$.
This defines a functor
\[
\Gamma_R:u\mathcal R\mathit{iesz}
\longrightarrow\mathcal{RMV},
\qquad (E,u)\longmapsto[0,u],
\]
which sends each morphism $f:(E,u)\to(F,v)$
to $\Gamma_R(f)=f|_{[0,u]}$.

\begin{theorem}[{\cite{DiNolaLeustean2014}}]
	\label{thm:riesz-mv-equivalence}
	The functor $\Gamma_R$ is an equivalence
	of categories.
\end{theorem}

In particular, we may choose a quasi-inverse
functor
$\Xi_R:\mathcal{RMV}\to u\mathcal R\mathit{iesz}$
whose underlying $\ell$-group construction
agrees with Mundici equivalence.

For a Riesz MV-algebra $A$, write
$\Xi_R(A)=(E_A,u_A)$.
The underlying unital $\ell$-group of $E_A$
is identified with the Mundici group
$(G_A,u_A)$ of the MV-reduct of $A$.
Under this identification, $\Xi_R(f)=\Xi(f)$
for every Riesz MV-homomorphism $f:A\to B$.
Moreover, the canonical MV-isomorphism
$\varphi_A:A\to[0,u_A]$
is also a Riesz MV-isomorphism, and these
isomorphisms form a natural isomorphism
$\operatorname{id}_{\mathcal{RMV}}\to\Gamma_R\Xi_R$.

A topology $\sigma$ on a real vector space $E$
is called a \emph{linear topology} if vector
addition and scalar multiplication are jointly
continuous, making $(E,\sigma)$ a topological
vector space.

A \emph{topological Riesz space} is a Riesz space
$E$ equipped with a linear topology $\sigma$
such that the lattice operations $\vee$ and
$\wedge$ are jointly continuous.
A subset $S$ of $E$ is called \emph{solid}
if $x\in S$ and $|y|\leq|x|$ imply $y\in S$.
A linear topology on $E$ is \emph{locally solid}
if it has a zero-neighbourhood basis consisting
of solid sets.
Much of the classical theory of topological
Riesz spaces is developed under the additional
assumption of local solidity
(see, e.g., \cite{AliprantisBurkinshaw2003}).
In this paper, we work in the general setting,
imposing local solidity only when needed.

A Riesz space $E$ is called \emph{Archimedean}
if $nx\leq y$ for every $n\geq1$ implies $x\leq0$.
If $E$ is Archimedean and has a strong unit $u$,
its \emph{order-unit norm} is defined by
\begin{equation}\label{eq:order-unit-norm}
	\|x\|_u
	=\inf\{\lambda>0:|x|\leq\lambda u\}.
\end{equation}
We denote the induced topology by $\nu_u$.
The closed unit ball is $[-u,u]$, and
the norm-bounded subsets of $E$ are precisely
its order-bounded subsets;
see \cite[Section~5.12]{NariciBeckenstein2011}.

Let $(V,\|\cdot\|)$ be a normed vector space
with norm topology $\nu$, and let
$\sigma\subseteq\nu$ be a Hausdorff
linear topology on $V$.
Wiweger~\cite{Wiweger1961} introduced the
\emph{mixed topology} $\gamma(\nu,\sigma)$
as the finest linear topology on $V$
satisfying
\[
\sigma\subseteq\gamma(\nu,\sigma)\subseteq\nu,
\qquad
\gamma(\nu,\sigma)|_B=\sigma|_B
\]
for every norm-bounded subset $B\subseteq V$.

Let $E$ be an Archimedean Riesz space with
strong unit $u$. Then $(E,\|\cdot\|_u)$
is a normed vector space.
Consequently, if $(E,\sigma)$ is a Hausdorff
topological Riesz space with
$\sigma\subseteq\nu_u$, the mixed topology
$\gamma(\nu_u,\sigma)$ is well defined.

In particular, every Hausdorff locally solid
topology $\sigma$ on $E$ satisfies
$\sigma\subseteq\nu_u$.
Indeed, every solid zero-neighbourhood
contains $[-\varepsilon u,\varepsilon u]$
for some $\varepsilon>0$.
Thus, the mixed topology
$\gamma(\nu_u,\sigma)$ is also well defined
in this case.

\begin{example}
	Let $X$ be a completely regular Hausdorff space.
	Then $C_b(X,\mathbb R)$, equipped with the
	pointwise order, is an Archimedean Riesz space
	with strong unit $\mathbf1$.
	Its order-unit norm is the supremum norm
	$\|\cdot\|_\infty$.
	
	The compact-open topology $\tau_c$ on
	$C_b(X,\mathbb R)$ is Hausdorff and locally
	convex-solid. The mixed topology
	$\gamma(\nu_{\mathbf1},\tau_c)$ coincides
	with the topology $\beta_0$, often called
	the strict topology;
	see \cite{Sentilles1972,FarkasKreidler2020}.
	
	The unit interval
	$[0,\mathbf1]=C_b(X,[0,1])$
	is a Riesz MV-algebra endowed with the
	induced compact-open topology.
\end{example}

We now fix the topological terminology used
throughout the paper. For a comprehensive
treatment of general topology, see
\cite{munkres2000}.

If $(X,\tau)$ is a topological space and
$Y\subseteq X$, we denote by $\tau|_Y$
the subspace topology induced by $\tau$.
A \emph{neighbourhood basis} at $x\in X$
is a family of neighbourhoods of $x$
such that every neighbourhood of $x$
contains some member of the family.
Its members need not be open.

We adopt the convention that the separation
properties regular ($T_3$), completely regular
($T_{3\frac12}$), and normal ($T_4$)
all include the $T_1$ condition.
Thus the separation axioms become
progressively stronger:
$T_4\Longrightarrow T_{3\frac12}
\Longrightarrow T_3\Longrightarrow T_2
\Longrightarrow T_1\Longrightarrow T_0.$

A family of subsets of a topological space
is \emph{locally finite} if every point has
a neighbourhood meeting only finitely many
members of the family.
A space is \emph{paracompact} if every open
cover admits a locally finite open refinement.
No separation axiom is included in this
definition.

We shall use the following classical result.

\begin{proposition}[{\cite[Theorem~41.1]{munkres2000}}]
	\label{prop:paracompact-hausdorff-normal}
	Every paracompact Hausdorff space is normal.
\end{proposition}

A space is \emph{$\sigma$-compact} if it is
a countable union of compact subspaces,
and \emph{pseudocompact} if every continuous
real-valued function on it is bounded.
It is \emph{$\sigma$-pseudocompact}
if it is a countable union of pseudocompact
subspaces, and \emph{locally $\sigma$-compact}
if every point has a $\sigma$-compact
neighbourhood.

A space has the \emph{countable chain condition}
(\emph{ccc}) if every pairwise disjoint
family of nonempty open sets is countable.
It has the \emph{Knaster property}
(or \emph{property K}) if every uncountable
family of nonempty open sets contains
an uncountable pairwise-intersecting subfamily.

\begin{proposition}[{\cite{FengGartside2026}}]
	\label{prop:knaster-implies-ccc}
	Every space with the Knaster property
	satisfies ccc.
\end{proposition}

Following \cite{FengGartside2026},
a family $\mathcal N$ of subsets of $X$
is a \emph{network for a cover} $\mathcal K$
if, whenever $K\in\mathcal K$ and $U$
is open with $K\subseteq U$,
there exists $N\in\mathcal N$ such that
$K\subseteq N\subseteq U$.
A space is a $\Sigma^{-}(\aleph_0)$-space
if it admits a cover by pseudocompact
subspaces having a countable network
for that cover.
In particular, every $\sigma$-pseudocompact
space is a $\Sigma^{-}(\aleph_0)$-space.

We next recall the basic notions concerning
topological groups.
A \emph{topological group} is a group
$(G,\cdot)$ equipped with a topology
such that multiplication and inversion
are continuous.
A \emph{group topology} on $G$ is a
topology making $G$ a topological group.
Every $T_0$ topological group is Hausdorff
and completely regular.

\begin{proposition}[{\cite[Theorem~3.1.1]{ArhangelskiiTkachenko2008}}]
	\label{prop:locally-sigma-compact-group}
	Every locally $\sigma$-compact Hausdorff
	topological group is paracompact.
\end{proposition}

In this paper, a \emph{lattice-ordered
	topological group} is an $\ell$-group
equipped with a group topology, without
any assumption that the lattice operations
are continuous.
We distinguish this notion from a
\emph{topological $\ell$-group}, in which
the lattice operations
\[
(x,y)\longmapsto x\vee y,
\qquad
(x,y)\longmapsto x\wedge y
\]
are also jointly continuous
(see, e.g., \cite{Madell1969,Ball1979,Jordan2023}).

We also recall the notion of retral spaces.
A continuous map $r:X\to Y$, where $Y$
is a subspace of $X$, is a \emph{retraction}
if $r(y)=y$ for every $y\in Y$.
A space is \emph{retral} if it is
homeomorphic to a retract of a
topological group
(see, e.g., \cite{GartsideReznichenkoSipacheva1997}).

The following result of Feng and Gartside
connects retral spaces with covering properties.

\begin{theorem}[{\cite[Theorem~4.11]{FengGartside2026}}]
	\label{thm:feng-gartside-knaster}
	Every $\Sigma^{-}(\aleph_0)$ retral space
	has the Knaster property.
\end{theorem}

Finally, recall the categorical formulation of
adjunctions. For functors
$F:\mathcal C\to\mathcal D$ and
$G:\mathcal D\to\mathcal C$, an \emph{adjunction}
$F\dashv G$ consists of natural bijections
\[
\operatorname{Hom}_{\mathcal D}(FX,Y)
\cong\operatorname{Hom}_{\mathcal C}(X,GY).
\]
Equivalently, it is specified by a unit
$\eta:\operatorname{id}_{\mathcal C}\to GF$
and a counit
$\varepsilon:FG\to\operatorname{id}_{\mathcal D}$
satisfying the triangle identities.
The left adjoint $F$ is fully faithful
if and only if the unit $\eta$ is a natural
isomorphism; see \cite{riehl2016}.

We conclude this section by fixing some notation
that will be used throughout the paper.
For an abelian $\ell$-group $G$ and $U\subseteq G$,
define
\begin{equation}\label{eq:U_triangle}
	U^\triangle
	:=\{g\in G:g^+\in U,\ g^-\in U\}.
\end{equation}
For a sequence $(U_i)_{i\geq1}$ of subsets
of an abelian group, define
\begin{equation}\label{eq:sum_Ui}
	\sum_{i=1}^{\infty}U_i
	:=\left\{
	\sum_{k=1}^n u_{i_k}:
	n\geq1,\ i_1<i_2<\cdots<i_n,\
	u_{i_k}\in U_{i_k}
	\right\}.
\end{equation}
Thus, despite the notation, only finite sums
are involved, with one element chosen from
each of finitely many distinct sets $U_i$.

\section{Topological lifting of Mundici equivalence}\label{sec:Topological_lifting}
Mundici equivalence represents every MV-algebra
as the unit interval of an abelian $\ell$-group
with a strong unit. When an MV-algebra carries
a topology, a natural question is whether this
representation can be made compatible with
the given topology.

In our previous work~\cite{JiaYang2026},
we studied the topological lifting problem
through the notion of summability and established
a corresponding lifting criterion.
We first recall the relevant definitions and
results before developing the canonical lifting
in this section.

\begin{definition}[{\cite{hoo1997}}]
	A \emph{topological MV-algebra} is an MV-algebra
	$(A,\oplus,\neg,0)$ equipped with a topology $\tau$
	such that $\oplus:A\times A\to A$ and
	$\neg:A\to A$ are continuous.
	No separation axiom is imposed on $\tau$.
\end{definition}

Since the derived MV-operations are term operations,
they are continuous on every topological MV-algebra.
In particular, the lattice operations $\vee$ and
$\wedge$, and the truncated difference $\ominus$,
are continuous.
We recall the following classical separation
properties of topological MV-algebras.

\begin{proposition}[{\cite{gan2024,hoo1997}}]
	\label{prop:basic-separation}
	Every $T_0$ topological MV-algebra is $T_3$.
\end{proposition}

We denote by
$$\mathcal T\!\mathit{op}\mathcal{MV}$$
the category of topological MV-algebras and
continuous MV-homomorphisms.

A fundamental feature of topological MV-algebras
is that their topology can be described using
neighbourhoods of $0$. For $U\subseteq A$ and $a\in A$,
put
\begin{equation}\label{eq:MV-neighbourhoods}
	U(a):=\{x\in A:x\ominus a\in U,\;
	a\ominus x\in U\}.
\end{equation}
If $\mathcal V(0)$ is a neighbourhood basis of $0$,
then
\[
\{U(a):U\in\mathcal V(0)\}
\]
is a neighbourhood basis of $a$ for every $a\in A$;
see~\cite{JiaYang2026}.
Here and below, a neighbourhood basis need not
consist of open sets.

\begin{definition}[{\cite{JiaYang2026}}]
	\label{def:summable}
	Let $(A,\tau)$ be a topological MV-algebra,
	and let $(G_A,u_A)=\Xi(A)$ be its Mundici group
	with canonical MV-isomorphism
	$\varphi_A:A\to[0,u_A]$
	as in Theorem~\ref{thm:mundici}.
	Let $\mathcal V(0)$ be a neighbourhood basis
	of $0$. The algebra $(A,\tau)$ is called
	\emph{summable} if, for every $U\in\mathcal V(0)$
	and $a\in A$, there is a sequence
	$(V_i)_{i\geq1}$ in $\mathcal V(0)$ such that
	\begin{equation}\tag{SC}\label{eq:SC}
		[0,u_A]\cap\left(
		\varphi_A(a)+\sum_{i=1}^{\infty}
		(\varphi_A(V_i))^\triangle
		\right)
		\subseteq\varphi_A(U(a)),
	\end{equation}
	where the notation is as in
	\eqref{eq:MV-neighbourhoods},
	\eqref{eq:U_triangle}, and
	\eqref{eq:sum_Ui}.
\end{definition}

Condition~\eqref{eq:SC}
does not depend on the choice of the
neighbourhood basis $\mathcal V(0)$.
It is an intrinsic condition on the topology
of $A$, although it is expressed using the
algebraic Mundici representation.

For a summable topological MV-algebra $(A,\tau)$,
our previous construction equips $G_A$ with the
group topology $\tau_{G_A}$ whose neighbourhood
basis at $0$ consists of the sets
\begin{equation}\label{eq:canonical-group-basis}
	\mathcal V_{G_A}
	:=\left\{
	\sum_{i=1}^{\infty}(\varphi_A(U_i))^\triangle
	:\ U_i\in\mathcal V(0)\ \text{for all }i\geq1
	\right\}.
\end{equation}

Strong summability provides the corresponding
criterion when the ambient group is required
to be a topological $\ell$-group
\cite[Proposition~4.9]{JiaYang2026}.

\begin{theorem}
	\label{thm:old-lifting}
	For a topological MV-algebra $(A,\tau)$,
	the following are equivalent:
	\begin{enumerate}
		\item $(A,\tau)$ is summable;
		\item there exist an abelian $\ell$-group $G$
		with a positive strong unit $u$ and a group
		topology $\sigma$ on $G$ such that $(A,\tau)$
		is MV-isomorphic and homeomorphic to
		\[
		\bigl([0,u],\sigma|_{[0,u]}\bigr).
		\]
	\end{enumerate}
\end{theorem}

The distinction between a lattice-ordered
topological group and a topological $\ell$-group
is important here. The preceding theorem requires
continuity of the group operations only, not of
the lattice operations. This leads to the stronger
notion introduced in~\cite{JiaYang2026}.

\begin{definition}[\cite{JiaYang2026}]
	\label{def:strongly-summable}
	A summable topological MV-algebra $(A,\tau)$
	is called \emph{strongly summable} if the
	group topology $\tau_{G_A}$ defined by
	\eqref{eq:canonical-group-basis} makes $G_A$
	a topological $\ell$-group.
\end{definition}

Summability provides a precise criterion for
representing a topological MV-algebra as the
unit interval of a lattice-ordered topological
group \cite[Theorem~3.5]{JiaYang2026}.

\begin{theorem}
	\label{thm:strong-lifting-criterion}
	For a summable topological MV-algebra
	$(A,\tau)$, the following are equivalent:
	\begin{enumerate}
		\item $(A,\tau)$ is strongly summable;
		\item there exist an abelian topological
		$\ell$-group $G$ with a positive strong unit
		$u$ such that $(A,\tau)$ is MV-isomorphic
		and homeomorphic to the unit interval
		$[0,u]$ with its subspace topology.
	\end{enumerate}
\end{theorem}

The preceding results yield a categorical
topological lifting of Mundici equivalence
for strongly summable topological MV-algebras,
in the form of an adjunction.

Let
\[
\mathcal T\!\mathit{op}^{\mathrm{str}}\!\mathcal{MV}
\]
denote the full subcategory of
$\mathcal T\!\mathit{op}\mathcal{MV}$
consisting of strongly summable objects,
and let
\[
\mathcal T\!\mathit{op}\,u\mathcal{LG}
\]
denote the category of abelian topological
$\ell$-groups with distinguished positive
strong unit and continuous unit-preserving
$\ell$-group homomorphisms.

The constructions $\Xi$ and $\Gamma$
define functors
\begin{equation}\label{eq:topological-functors}
	\Xi:
	\mathcal T\!\mathit{op}^{\mathrm{str}}\!\mathcal{MV}
	\rightleftarrows
	\mathcal T\!\mathit{op}\,u\mathcal{LG}
	:\Gamma.
\end{equation}
On objects, they are given by
\[
\Xi(A,\tau)=(G_A,u_A,\tau_{G_A}),
\qquad
\Gamma(G,u,\sigma)
=\bigl([0,u],\sigma|_{[0,u]}\bigr).
\]
On morphisms, $\Xi$ acts by the canonical
extension of MV-homomorphisms, while
$\Gamma$ acts by restriction to unit intervals.

We recall the following result from
\cite[Corollary~4.13]{JiaYang2026}.

\begin{proposition}[{\cite[Corollary~4.13]{JiaYang2026}}]
	\label{prop:strong-adjunction}
	The functors in \eqref{eq:topological-functors}
	form an adjunction
	\[
	\Xi\dashv\Gamma.
	\]
	Its unit has components
	\[
	\varphi_A:(A,\tau)\longrightarrow
	\bigl([0,u_A],\tau_{G_A}|_{[0,u_A]}\bigr),
	\]
	which are MV-isomorphisms and homeomorphisms.
	Consequently, $\Xi$ is fully faithful.
\end{proposition}

This adjunction is not, in general, an
equivalence of categories;
see \cite[Example~4.14]{JiaYang2026}.

We now extend the preceding topological lifting
to all topological MV-algebras by proving that
every topological MV-algebra is strongly summable.
In view of Theorems~\ref{thm:old-lifting} and
\ref{thm:strong-lifting-criterion}, we proceed
in two steps: first establishing summability
and then strong summability.
The following elementary decomposition
will be useful in both steps.

\begin{lemma}
	\label{lem:layer-decomposition}
	Let $(G,u)$ be an abelian $\ell$-group with strong unit. For
	$x\in G$ and $k\geq1$, define
	\begin{equation}
		\label{eq:layer-decomposition-maps}
		P_k(x):=\bigl((x^+-(k-1)u)^+\bigr)\wedge u,
		\qquad
		N_k(x):=\bigl((x^--(k-1)u)^+\bigr)\wedge u.
	\end{equation}
	If $|x|\leq nu$ for some $n\geq1$, then
	\[
	x^+=\sum_{k=1}^n P_k(x),\qquad
	x^-=\sum_{k=1}^n N_k(x),
	\]
	and consequently,
	\[
	x=\sum_{k=1}^n P_k(x)-\sum_{k=1}^n N_k(x),
	\qquad P_k(x),N_k(x)\in[0,u].
	\]
\end{lemma}

\begin{proof}
	For any $a\in G$, we have
	\begin{equation}\label{eq:truncation-identity}
		a\wedge u=a+(0\wedge(u-a))
		=a-(0\vee(a-u))=a-(a-u)^+.
	\end{equation}
	Moreover, since $u\geq0$,
	\begin{equation}\label{eq:shifted-positive-part}
		(a^+-u)^+=((a-u)\vee(-u))\vee0=(a-u)^+.
	\end{equation}
	Now let $0\leq a\leq nu$ and put $a_k=(a-(k-1)u)^+$ for $1\leq k\leq n+1$. By \eqref{eq:shifted-positive-part}, $(a_k-u)^+=a_{k+1}$, so \eqref{eq:truncation-identity} gives $a_k=(a_k\wedge u)+a_{k+1}$. Applying this recursively and using $a_1=a$ and $a_{n+1}=0$, we obtain
	\[
	a=\sum_{k=1}^n(a_k\wedge u).
	\]
	Since $0\leq x^+,x^-\leq nu$, applying this formula to $a=x^+$ and $a=x^-$ gives the asserted decompositions of $x^+$ and $x^-$. The decomposition of $x$ follows from $x=x^+-x^-$, while $P_k(x),N_k(x)\in[0,u]$ holds by definition.
\end{proof}

\begin{proposition}
	\label{prop:universal-summability}
	Every topological MV-algebra is summable.
\end{proposition}

\begin{proof}
	Let $(A,\tau)$ be a topological MV-algebra, and let
	$(G,u)=\Xi(A)$ be its Mundici group. By Mundici equivalence,
	we may assume without loss of generality that
	$A=\Gamma(G,u)=[0,u]$.
	Let $\mathcal V(0)$ denote the family of zero-neighbourhoods
	in $(A,\tau)$.
	
	Under this identification, the summability condition
	\eqref{eq:SC} reduces to the following: for every
	$U\in\mathcal V(0)$ and $a\in A$, there exist
	$V_i\in\mathcal V(0)$ such that
	\[
	[0,u]\cap
	\left(a+\sum_{i=1}^{\infty}V_i^\triangle\right)
	\subseteq U(a).
	\]
	
	We shall prove the following stronger claim.
	
	\textbf{Claim:}\quad
	For every $U\in\mathcal V(0)$, there exist
	$V_i\in\mathcal V(0)$ such that
	\[
	t^+\wedge u,\ t^-\wedge u\in U
	\qquad
	\text{for all }t\in\sum_{i=1}^{\infty}V_i^\triangle.
	\]
	
	Indeed, if the claim holds, then for any $a\in A$ and
	$x\in[0,u]\cap(a+\sum_{i=1}^{\infty}V_i^\triangle)$,
	putting $t=x-a$ gives $t^+,t^-\leq u$. Hence
	$x\ominus a=t^+\wedge u\in U$ and
	$a\ominus x=t^-\wedge u\in U$.
	Thus $x\in U(a)$, which verifies \eqref{eq:SC}.

	We now prove the claim. For $n\geq1$ and $t\in G$,
	using the maps $P_k$ and $N_k$ from
	\eqref{eq:layer-decomposition-maps}, put
	\begin{equation}
		\label{eq:sigma-definition}
		\begin{aligned}
			c_n(t)&:=(-nu)\vee(t\wedge nu),\\
			P_k(t)&:=\bigl((t^+-(k-1)u)^+\bigr)\wedge u,\\
			N_k(t)&:=\bigl((t^--(k-1)u)^+\bigr)\wedge u,\\
			\sigma_n(t)&:=
			\bigl(P_1(t),P_2(t),\ldots,P_n(t),
			N_1(t),N_2(t),\ldots,N_n(t)\bigr)\in A^{2n}.
		\end{aligned}
	\end{equation}
	
	For every $t\in G$ and $h\in G$ with $|h|\leq u$,
	we have
	\[
	c_n(h+t)=c_n\!\left(h+c_{n+1}(t)\right),
	\qquad
	\sigma_n(c_n(t))=\sigma_n(t).
	\]
	By Proposition~\ref{prop:subdirect-representation}(4), it suffices to
	verify these identities in a totally ordered $\ell$-group.
	The first follows by a straightforward calculation in the
	three cases $t\geq(n+1)u$, $|t|\leq(n+1)u$, and
	$t\leq-(n+1)u$, while the second follows directly from
	the definitions of $P_k$ and $N_k$.
	Consequently,
	\[
	\begin{aligned}
		\sigma_n(h+t)
		&=\sigma_n(c_n(h+t))\\
		&=\sigma_n\!\left(c_n(h+c_{n+1}(t))\right)\\
		&=\sigma_n\!\left(h+c_{n+1}(t)\right).
	\end{aligned}
	\]
	Moreover, applying Lemma~\ref{lem:layer-decomposition}
	to $c_{n+1}(t)$ and using
	$\sigma_{n+1}(c_{n+1}(t))=\sigma_{n+1}(t)$,
	we obtain
	\[
	c_{n+1}(t)
	=
	\sum_{k=1}^{n+1}P_k(t)
	-
	\sum_{k=1}^{n+1}N_k(t).
	\]
	
	For $p,q\in A$ and
	$\mathbf r=(r_1,r_2,\ldots,r_{n+1}),
	\mathbf s=(s_1,s_2,\ldots,s_{n+1})\in A^{n+1}$,
	put $d=p-q+\sum_{j=1}^{n+1}(r_j-s_j).$
	Define
	\begin{equation}
		\label{eq:Lambda-definition}
		\begin{aligned}
			\Lambda_n:
			A^2\times A^{2(n+1)}
			&\longrightarrow A^{2n},\\
			(p,q,\mathbf r,\mathbf s)
			&\longmapsto
			\bigl((P_k(d))_{k=1}^n,(N_k(d))_{k=1}^n\bigr)\\
			&=
			\left(
			\bigl((d^+-(k-1)u)^+\wedge u\bigr)_{k=1}^n,
			\bigl((d^--(k-1)u)^+\wedge u\bigr)_{k=1}^n
			\right).
		\end{aligned}
	\end{equation}
	To simplify the coordinate expressions of $\Lambda_n$,
	we observe that
	\[
	(z^+-a)^+=(z-a)^+
	\qquad (z\in G,\ a\in G^+).
	\]
	Indeed, by Proposition~\ref{prop:subdirect-representation}(4),
	it suffices to verify this identity in a totally
	ordered $\ell$-group. If $z\geq0$, the identity is
	immediate. If $z\leq0$, both sides equal $0$.
	
	Taking $a=(k-1)u$ and applying this identity to
	$z=d$ and $z=-d$, respectively, we obtain
	\[
	\begin{aligned}
		P_k(d)
		&=\bigl((d^+-(k-1)u)^+\bigr)\wedge u\\
		&=(d-(k-1)u)^+\wedge u,\\
		N_k(d)
		&=\bigl((d^--(k-1)u)^+\bigr)\wedge u\\
		&=\bigl(((-d)^+-(k-1)u)^+\bigr)\wedge u\\
		&=(-d-(k-1)u)^+\wedge u.
	\end{aligned}
	\]
	Consequently,
	\[
	\Lambda_n(p,q,\mathbf r,\mathbf s)
	=
	\left(
	\bigl((d-(k-1)u)^+\wedge u\bigr)_{k=1}^n,
	\bigl((-d-(k-1)u)^+\wedge u\bigr)_{k=1}^n
	\right).
	\]
	
	By Lemma~\ref{lem:integer-affine-mv-term}, each coordinate
	of $\Lambda_n$ is an MV-term function. Hence
	\[
	\Lambda_n:
	(A^2\times A^{2(n+1)},\tau^{2(n+2)})
	\longrightarrow
	(A^{2n},\tau^{2n})
	\]
	is continuous. 
	Moreover, $\Lambda_n(0)=0$ by
	\eqref{eq:Lambda-definition}.
	Using the decomposition of $c_{n+1}(t)$
	established above, we obtain, for $|h|\leq u$,
	\begin{equation}
		\label{eq:Lambda-identity}
		\begin{aligned}
			\Lambda_n(h^+,h^-,\sigma_{n+1}(t))
			&=\sigma_n(h+c_{n+1}(t))\\
			&=\sigma_n(h+t).
		\end{aligned}
	\end{equation}
	
	Fix $U\in\mathcal V(0)$ and put $O_1=U\times U$.
	Since $\Lambda_1$ is continuous at $0$ and
	$\Lambda_1(0)=0$, we can choose $V_1\in\mathcal V(0)$
	and a product zero-neighbourhood $O_2\subseteq A^4$
	such that
	$\Lambda_1(V_1\times V_1\times O_2)\subseteq O_1$.
	
	Inductively, suppose that a product zero-neighbourhood
	$O_n\subseteq A^{2n}$ has been chosen.
	By the continuity of $\Lambda_n$ at $0$ and
	$\Lambda_n(0)=0$, we can similarly choose
	$V_n\in\mathcal V(0)$ and a product zero-neighbourhood
	$O_{n+1}\subseteq A^{2(n+1)}$ such that
	\begin{equation}
		\label{eq:Lambda-neighborhood}
		\Lambda_n(V_n\times V_n\times O_{n+1})
		\subseteq O_n
		\qquad(n\geq1).
	\end{equation}
	
	Take any $t\in\sum_{i=1}^{\infty}V_i^\triangle$.
	By definition, we may write $t=\sum_{i=1}^N h_i$,
	where $h_i^+,h_i^-\in V_i$, inserting zero terms
	at omitted indices.
	Put $t_i=\sum_{j=i}^N h_j$ and $t_{N+1}=0$.
	Since $h_i^+,h_i^-\in A$, we have $|h_i|\leq u$.
	
	We show that $\sigma_i(t_i)\in O_i$ by backward induction.
	The initial case follows from
	$\sigma_{N+1}(t_{N+1})=0\in O_{N+1}$.
	Now suppose that $\sigma_{i+1}(t_{i+1})\in O_{i+1}$.
	Since $t_i=h_i+t_{i+1}$ and $h_i^+,h_i^-\in V_i$,
	by \eqref{eq:Lambda-identity} and
	\eqref{eq:Lambda-neighborhood}, we obtain
	\[
	\sigma_i(t_i)
	=\Lambda_i(h_i^+,h_i^-,\sigma_{i+1}(t_{i+1}))
	\in O_i.
	\]
	Thus the induction is complete. In particular,
	\[
	\sigma_1(t)
	=(t^+\wedge u,t^-\wedge u)
	\in O_1=U\times U.
	\]
	This proves the claim and hence verifies
	\eqref{eq:SC}. Therefore $(A,\tau)$ is summable.
\end{proof}

We next strengthen this result by showing that every
topological MV-algebra is strongly summable. The key
step is to establish the continuity of the lattice
operations on its Mundici group.

\begin{proposition}
	\label{prop:universal-strong-summability}
	Every topological MV-algebra is strongly summable.
\end{proposition}

\begin{proof}
	Let $(A,\tau)$ be a topological MV-algebra, and let
	$(G,u)=\Xi(A)$ be its Mundici group. By Mundici equivalence,
	we may assume without loss of generality that
	$A=\Gamma(G,u)=[0,u]$.
	By Proposition~\ref{prop:universal-summability},
	$(A,\tau)$ is summable. Let $\tau_G$ be the induced
	group topology on $G$, and let $\mathcal V(0)$ denote
	the family of zero-neighbourhoods in $(A,\tau)$.
	By construction, the sets
	$W=\sum_{i=1}^{\infty}U_i^\triangle$,
	where $U_i\in\mathcal V(0)$, form a neighbourhood
	basis at $0$ in $(G,\tau_G)$.
	
	Since $g\vee h=h+(g-h)^+$ and
	$g\wedge h=g-(g-h)^+$, it suffices to prove
	that the positive-part map $g\mapsto g^+$
	is continuous on $(G,\tau_G)$.
	
	To this end, fix $x\in G$ and a sequence
	$(U_i)_{i\geq1}$ in $\mathcal V(0)$.
	To control changes in the positive-part operation,
	define the increment map
	\begin{equation}
		\label{eq:delta-definition}
		\begin{aligned}
			\delta:G\times G&\longrightarrow G,\\
			(y,h)&\longmapsto (y+h)^+-y^+.
		\end{aligned}
	\end{equation}
	For any $N\geq1$ and $h_1,h_2,\ldots,h_N\in G$,
	put $t_i=\sum_{j=i}^N h_j$ and $t_{N+1}=0$.
	Then
	\[
	\left(x+\sum_{i=1}^N h_i\right)^+-x^+
	=\sum_{i=1}^N\delta(x+t_{i+1},h_i).
	\]
	Thus we shall construct a sequence $(V_i)_{i\geq1}$
	in $\mathcal V(0)$ such that, whenever
	$h_i^+,h_i^-\in V_i$ for $1\leq i\leq N$,
	each increment $\delta(x+t_{i+1},h_i)$
	belongs to $U_i^\triangle$.
	This will establish the required inclusion
	\[
	\left(x+\sum_{i=1}^{\infty}V_i^\triangle\right)^+
	\subseteq
	x^++\sum_{i=1}^{\infty}U_i^\triangle.
	\]
	
	We retain the maps $\sigma_n$ and $\Lambda_n$ defined
	in \eqref{eq:sigma-definition} and
	\eqref{eq:Lambda-definition}, respectively.
	Recall that $\Lambda_n$ is continuous and satisfies
	\eqref{eq:Lambda-identity}, namely,
	\[
	\Lambda_n(h^+,h^-,\sigma_{n+1}(t))
	=\sigma_n(h+t)
	\qquad(t\in G,\ |h|\leq u).
	\]
	
	To express the positive and negative parts
	of the increments defined in
	\eqref{eq:delta-definition}, introduce
	the following auxiliary map:
	\begin{equation}
		\label{eq:D-definition}
		\begin{aligned}
			D:A^4&\longrightarrow A^2,\\
			(a,b,p,q)&\longmapsto
			\bigl(a\ominus(b\oplus q),\,b\wedge p\bigr).
		\end{aligned}
	\end{equation}
	For $y,h\in G$ with $|h|\leq u$, we claim that
	\begin{equation}
		\label{eq:D-increment}
		\bigl(\delta(y,h)^+,\delta(y,h)^-\bigr)
		=D(h^+,h^-,y^+\wedge u,y^-\wedge u).
	\end{equation}
	By Proposition~\ref{prop:subdirect-representation}(4),
	it suffices to verify \eqref{eq:D-increment}
	in a totally ordered $\ell$-group.
	
	If $h\geq0$, then $0\leq h\leq u$.
	When $y\geq0$, we have $\delta(y,h)=h$;
	when $y\leq0$, we have
	$\delta(y,h)=(h-y^-)^+$.
	Thus, in both cases,
	$\delta(y,h)=(h-y^-)^+
	=(h-(y^-\wedge u))^+
	=h\ominus(y^-\wedge u)$,
	where the second equality follows from $h\leq u$.
	Hence $\delta(y,h)^+=h\ominus(y^-\wedge u)$
	and $\delta(y,h)^-=0$.
	
	If $h\leq0$, then $0\leq-h\leq u$.
	When $y\leq0$, we have $\delta(y,h)=0$;
	when $y\geq0$, we have
	$\delta(y,h)=(y+h)^+-y=-((-h)\wedge y)$.
	Thus, in both cases,
	$\delta(y,h)=-((-h)\wedge y^+)$.
	Hence $\delta(y,h)^+=0$ and
	$\delta(y,h)^-=(-h)\wedge y^+
	=(-h)\wedge(y^+\wedge u)$,
	where the last equality follows from $-h\leq u$.
	
	In either case, \eqref{eq:D-increment} follows
	from the definition of $D$ in \eqref{eq:D-definition}.
	Since $D$ is an MV-term map, it is continuous.
	Moreover, $D(0,0,p,q)=(0,0)$ for all $p,q\in A$.
	
	For the fixed $x\in G$, choose $m\geq1$
	such that $|x|\leq mu$.
	Using the maps $P_k$ and $N_k$ defined in
	\eqref{eq:layer-decomposition-maps},
	put $e_k=P_k(x)-N_k(x)$ for $1\leq k\leq m$.
	By Lemma~\ref{lem:layer-decomposition},
	we have $x=\sum_{k=1}^m e_k$.
	Since $P_k(x)\wedge N_k(x)=0$,
	we also have $e_k^+=P_k(x)$,
	$e_k^-=N_k(x)$, and $|e_k|\leq u$.
	
	We now construct a continuous map
	$F_x:A^{2(m+1)}\to A^2$ by backward recursion.
	Set $F_{m+1}=\mathrm{id}_{A^{2(m+1)}}$.
	For $k=m,m-1,\ldots,1$, define
	$F_k:A^{2(m+1)}\to A^{2k}$ by
	\[
	F_k(z):=
	\Lambda_k(e_k^+,e_k^-,F_{k+1}(z)),
	\qquad z\in A^{2(m+1)}.
	\]
	Since each $\Lambda_k$ is continuous,
	all the maps $F_k$ are continuous.
	
	We show by backward induction that
	$F_k(\sigma_{m+1}(t))
	=\sigma_k(t+\sum_{j=k}^m e_j)$
	for every $t\in G$ and $1\leq k\leq m+1$,
	where the sum is understood to be zero when $k=m+1$.
	The initial case $k=m+1$ follows from
	the definition of $F_{m+1}$.
	Now suppose the identity holds for $k+1$.
	Since $|e_k|\leq u$, by
	\eqref{eq:Lambda-identity}, we obtain
	\[
	\begin{aligned}
		F_k(\sigma_{m+1}(t))
		&=\Lambda_k\left(
		e_k^+,e_k^-,
		\sigma_{k+1}\left(t+\sum_{j=k+1}^m e_j\right)
		\right)\\
		&=\sigma_k\left(t+\sum_{j=k}^m e_j\right).
	\end{aligned}
	\]
	Thus the induction is complete.
	
	Taking $F_x=F_1$ and using
	$x=\sum_{k=1}^m e_k$, we obtain
	\begin{equation}
		\label{eq:F-x-identity}
		F_x(\sigma_{m+1}(t))
		=\sigma_1(x+t)
		=\bigl((x+t)^+\wedge u,(x+t)^-\wedge u\bigr),
	\end{equation}
	for every $t\in G$.
	
	We next extend this construction to higher-dimensional
	coordinate spaces by retaining only the first $m+1$
	coordinates in each of the two blocks.
	For $n\geq m+1$ and
	$\mathbf v=(v_1,v_2,\ldots,v_n),
	\mathbf w=(w_1,w_2,\ldots,w_n)\in A^n$,
	define
	\begin{equation}
		\label{eq:Delta-definition}
		\begin{aligned}
			\Delta_{x,n}:A^2\times A^{2n}
			&\longrightarrow A^2,\\
			(a,b,\mathbf v,\mathbf w)
			&\longmapsto
			D\bigl(a,b,
			F_x(v_1,v_2,\ldots,v_{m+1},
			w_1,w_2,\ldots,w_{m+1})\bigr).
		\end{aligned}
	\end{equation}
	Since coordinate projections, $F_x$, and $D$
	are continuous, so is $\Delta_{x,n}$.
	Moreover, $\Delta_{x,n}(0)=0$, since
	$D(0,0,p,q)=(0,0)$ for all $p,q\in A$.
	
	For $t,h\in G$ with $|h|\leq u$,
	by \eqref{eq:sigma-definition},
	\eqref{eq:D-increment},
	\eqref{eq:F-x-identity}, and
	\eqref{eq:Delta-definition}, we obtain
	\begin{equation}
		\label{eq:Delta-increment}
		\Delta_{x,n}(h^+,h^-,\sigma_n(t))
		=
		\bigl(\delta(x+t,h)^+,\delta(x+t,h)^-\bigr).
	\end{equation}
	
	We now construct the sequences $(V_i)_{i\geq1}$
	and $(O_i)_{i\geq1}$ recursively.
	Set $O_1=A^{2(m+1)}$.
	Since the joint map
	$(\Lambda_{m+1},\Delta_{x,m+2})$ is continuous
	at $0$ and maps $0$ to $0$, we can choose
	$V_1\in\mathcal V(0)$ and a product
	zero-neighbourhood $O_2\subseteq A^{2(m+2)}$
	such that the image of
	$V_1\times V_1\times O_2$ is contained in
	$O_1\times U_1\times U_1$.
	
	Inductively, suppose that a product
	zero-neighbourhood
	$O_i\subseteq A^{2(m+i)}$ has been chosen.
	The joint map
	\[
	(\Lambda_{m+i},\Delta_{x,m+i+1}):
	A^2\times A^{2(m+i+1)}
	\longrightarrow A^{2(m+i)}\times A^2
	\]
	is continuous at $0$ and maps $0$ to $0$.
	Hence we can similarly choose
	$V_i\in\mathcal V(0)$ and a product
	zero-neighbourhood
	$O_{i+1}\subseteq A^{2(m+i+1)}$ such that
	\begin{equation}
		\label{eq:Lambda-Delta-neighborhood}
		(\Lambda_{m+i},\Delta_{x,m+i+1})
		(V_i\times V_i\times O_{i+1})
		\subseteq O_i\times U_i\times U_i
		\qquad(i\geq1).
	\end{equation}
	
	Take any $h\in\sum_{i=1}^{\infty}V_i^\triangle$.
	By definition, we may write $h=\sum_{i=1}^N h_i$,
	where $h_i^+,h_i^-\in V_i$, inserting zero terms
	at omitted indices.
	Put $t_i=\sum_{j=i}^N h_j$ and $t_{N+1}=0$.
	Since $h_i^+,h_i^-\in A$, we have $|h_i|\leq u$.
	
	We show that $\sigma_{m+i}(t_i)\in O_i$
	by backward induction.
	The initial case follows from
	$\sigma_{m+N+1}(t_{N+1})=0\in O_{N+1}$.
	Now suppose that
	$\sigma_{m+i+1}(t_{i+1})\in O_{i+1}$.
	Since $t_i=h_i+t_{i+1}$ and
	$h_i^+,h_i^-\in V_i$,
	by \eqref{eq:Lambda-identity} and the first
	component of \eqref{eq:Lambda-Delta-neighborhood},
	we obtain
	\[
	\sigma_{m+i}(t_i)
	=\Lambda_{m+i}
	(h_i^+,h_i^-,\sigma_{m+i+1}(t_{i+1}))
	\in O_i.
	\]
	Thus the induction is complete.

	For each $1\leq i\leq N$, we have
	$\sigma_{m+i+1}(t_{i+1})\in O_{i+1}$ by
	the backward induction above.
	Since $h_i^+,h_i^-\in V_i$,
	\eqref{eq:Delta-increment} and the second
	component of \eqref{eq:Lambda-Delta-neighborhood}
	give
	$\delta(x+t_{i+1},h_i)^+,
	\delta(x+t_{i+1},h_i)^-\in U_i$.
	Hence $\delta(x+t_{i+1},h_i)\in U_i^\triangle$.

	Since $t_1=h$ and $t_{N+1}=0$, we obtain
	\[
	(x+h)^+-x^+
	=\sum_{i=1}^N\delta(x+t_{i+1},h_i)
	\in\sum_{i=1}^{\infty}U_i^\triangle.
	\]
	This establishes the required inclusion
	and proves that the positive-part map is continuous
	on $(G,\tau_G)$.
	Consequently, both lattice operations are continuous,
	so $(G,\tau_G)$ is a topological $\ell$-group.
	Therefore $(A,\tau)$ is strongly summable.
\end{proof}

Combining Propositions~\ref{prop:universal-summability}
and~\ref{prop:universal-strong-summability}
with Theorems~\ref{thm:old-lifting}
and~\ref{thm:strong-lifting-criterion},
we obtain the main result of this paper.

\begin{theorem}
	\label{thm:topological-lifting}
	Let $(A,\tau)$ be a topological MV-algebra, and let
	$(G_A,u_A)=\Xi(A)$ be its Mundici group.
	Then there exists a group topology $\tau_{G_A}$ on
	$G_A$ such that $(G_A,\tau_{G_A})$ is a topological
	$\ell$-group and the canonical MV-isomorphism
	\[
	\varphi_A:(A,\tau)\longrightarrow
	([0,u_A],\tau_{G_A}|_{[0,u_A]})
	\]
	is a homeomorphism.
	
	In particular, every topological MV-algebra is
	MV-isomorphic and homeomorphic to the unit interval
	of a topological $\ell$-group.
\end{theorem}

We now establish the categorical counterpart of
Theorem~\ref{thm:topological-lifting}.

Recall that $\mathcal T\!\mathit{op}\mathcal{MV}$
denotes the category of topological MV-algebras
and continuous MV-homomorphisms, while
$\mathcal T\!\mathit{op}\,u\mathcal{LG}$
denotes the category of abelian topological
$\ell$-groups with strong unit and continuous
unit-preserving $\ell$-group homomorphisms.

For the functors $\Xi$ and $\Gamma$ introduced in
\eqref{eq:topological-functors}, combining
Propositions~\ref{prop:universal-strong-summability}
and~\ref{prop:strong-adjunction}
yields the following result.

\begin{theorem}[Topological lifting of Mundici equivalence]
	\label{thm:categorical-topological-lifting}
	The functors $\Xi$ and $\Gamma$ induce an adjunction
	\[
	\Xi\dashv\Gamma:
	\mathcal T\!\mathit{op}\mathcal{MV}
	\rightharpoonup
	\mathcal T\!\mathit{op}\,u\mathcal{LG}.
	\]
	Moreover, the unit $\eta_A:A\to\Gamma\Xi(A)$
	is a natural isomorphism. Consequently, $\Xi$
	is fully faithful, and
	$\mathcal T\!\mathit{op}\mathcal{MV}$
	is equivalent to the full subcategory of
	$\mathcal T\!\mathit{op}\,u\mathcal{LG}$
	given by the essential image of $\Xi$.
\end{theorem}

\begin{remark}\label{rem:canonical-topology}
	\noindent\textup{(1)}
	The canonical topology $\tau_{G_A}$ is the finest
	group topology on $G_A$ inducing $\tau$ on
	$A\cong[0,u_A]$.
	
	\smallskip
	\noindent\textup{(2)}
	Every topology induces a convergence, but not every
	convergence is topological. Ball, Georgescu and
	Leu\c{s}tean~\cite{BallGeorgescuLeustean2002}
	established a one-to-one correspondence between
	MV-convergences on $A=[0,u]$ and bounded
	$\ell_u$-convergences on its Mundici group $G$.
	However, this correspondence does not necessarily
	preserve topologicality. In particular, the
	topologicality equivalence stated in
	\cite[Remark~3.4(2)]{BallGeorgescuLeustean2002}
	does not hold in full generality, as the following
	example shows.
	
	Let $A=[0,1]^{\mathbb N}$ carry the product topology,
	and let $G=\ell^\infty(\mathbb N,\mathbb R)$ with
	$u=(1,1,\ldots)$. Write $\to_A$ for the convergence
	induced by the product topology on $A$, and
	$\to_G$ for its bounded extension to $G$ under
	the correspondence of
	\cite[Theorem~3.3]{BallGeorgescuLeustean2002}.
	
	For each $n\geq1$, put $B_n=[-nu,nu]$, and let
	$\mathcal P_n$ be the filter generated by the sets
	$B_n\cap U$, where $U$ ranges over the basic product
	neighbourhoods of $0$ in $\mathbb R^{\mathbb N}$.
	By the bounded-extension construction in
	\cite[Theorem~3.3]{BallGeorgescuLeustean2002},
	we have
	\[
	\mathcal F\to_G 0
	\quad\Longleftrightarrow\quad
	\mathcal F\supseteq\mathcal P_n
	\text{ for some }n\geq1.
	\]
	Moreover,
	\[
	\mathcal P_1\supsetneq\mathcal P_2
	\supsetneq\cdots.
	\]
	Indeed, $\mathcal P_n\supseteq\mathcal P_{n+1}$,
	and for every basic product neighbourhood $U$
	one can choose a coordinate $j$ not restricted
	by $U$, so that
	$(n+1)e_j\in(B_{n+1}\cap U)\setminus B_n$.
	Thus $B_n\in\mathcal P_n\setminus\mathcal P_{n+1}$.
	
	By Ball's topologicality criterion
	\cite[Section~1]{Ball1980},
	if $\to_G$ were topological, there would exist
	a coarsest filter $\mathcal N_G(0)$ converging to $0$.
	By the characterization above,
	$\mathcal N_G(0)\supseteq\mathcal P_m$
	for some $m$, while
	$\mathcal N_G(0)\subseteq\mathcal P_{m+1}$,
	since $\mathcal P_{m+1}\to_G0$.
	This contradicts
	$\mathcal P_m\supsetneq\mathcal P_{m+1}$.
	Therefore the bounded extension of a topological
	MV-convergence need not be topological.
\end{remark}

\section{Applications}\label{sec:Applications}

In this section, we explore several applications of the
topological lifting established in
Section~\ref{sec:Topological_lifting}.
We first investigate topological properties of
MV-algebras, including separation axioms
(Corollary~\ref{cor:separation-equivalence}),
and establish their canonical retract representation
(Proposition~\ref{prop:canonical-retract}),
which leads to further consequences for covering
properties under compactness-related assumptions.

We then extend the lifting to topological Riesz MV-algebras
(Proposition~\ref{prop:topological-riesz-lifting})
and establish a connection between the canonical lifted
topology and mixed topologies on Riesz spaces
(Proposition~\ref{prop:canonical-mixedification}).
Building on the lifting of strongly continuous semigroups
(Proposition~\ref{prop:mv-semigroup-lifting}),
we obtain a unit-interval Koopman representation theorem
(Theorem~\ref{thm:unit-interval-koopman}).

\subsection{Topological properties of MV-algebras}
\label{subsec:topological-properties}

Throughout this subsection, for a topological
MV-algebra $(A,\tau)$, we write
$\Xi(A,\tau)=(G_A,u_A,\tau_{G_A})$
for its canonical topological $\ell$-group
and identify $A$ with $[0,u_A]$ via $\varphi_A$.

We first establish two basic topological properties
preserved by the canonical lifting.

\begin{proposition}\label{prop:canonical-transfer}
	Let $(A,\tau)$ be a topological MV-algebra. Then:
	\begin{enumerate}
		\item $(A,\tau)$ is $T_0$ if and only if
		$(G_A,\tau_{G_A})$ is $T_0$, equivalently, Hausdorff;
		\item $(A,\tau)$ is connected if and only if
		$(G_A,\tau_{G_A})$ is connected.
	\end{enumerate}
\end{proposition}

\begin{proof}
	For (1), one implication follows from heredity of $T_0$.
	Conversely, suppose that $A$ is $T_0$.
	By Proposition~\ref{prop:basic-separation}, $A$ is $T_1$.
	Put $N=\overline{\{0\}}^{\,G_A}$.
	Since $A=[0,u_A]$ carries the subspace topology,
	we have $N\cap[0,u_A]=\{0\}$.
	
	If $0\neq g\in N$, continuity of the lattice
	operations gives $|g|\wedge u_A\in N$.
	Moreover, $|g|\wedge u_A>0$, since otherwise
	the strong-unit property together with the
	Riesz decomposition property would imply $g=0$.
	Thus $N\cap[0,u_A]\neq\{0\}$, a contradiction.
	Hence $N=\{0\}$, so $G_A$ is Hausdorff.
	
	For (2), if $A$ is connected, then
	$G_A=\bigcup_{n\geq1}n(A-A)$ by the strong-unit
	property and the Riesz decomposition property.
	Each $n(A-A)$ is a connected continuous image
	of $A^{2n}$, and these sets form an increasing
	family; hence $G_A$ is connected.
	Conversely, if $G_A$ is connected,
	the continuous surjection
	$r_A(g)=g^+\wedge u_A$ from $G_A$ onto
	$A$ shows that $A$ is connected.
\end{proof}

The separation properties of topological MV-algebras
have been studied since Hoo~\cite{hoo1997}, who proved
the equivalence of $T_0$, $T_1$, and $T_2$.
Gan et al.~\cite{gan2024} subsequently extended
this equivalence to $T_3$.
More recently, Li and Yang~\cite{li2026} raised the
question of whether every $T_0$ topological MV-algebra
is completely regular, and answered it affirmatively
under the assumption of local convexity.
The following corollary resolves this question
in full generality.

\begin{corollary}\label{cor:separation-equivalence}
	For a topological MV-algebra $(A,\tau)$, the following conditions are
	equivalent:
	\begin{enumerate}
		\item $(A,\tau)$ is $T_0$;
		\item $(A,\tau)$ is $T_1$;
		\item $(A,\tau)$ is Hausdorff;
		\item $(A,\tau)$ is regular;
		\item $(A,\tau)$ is completely regular.
	\end{enumerate}
\end{corollary}

\begin{proof}
	It suffices to prove that $T_0$ implies complete regularity.
	By Proposition~\ref{prop:canonical-transfer}, $G_A$ is a Hausdorff
	topological group and hence Tychonoff. Therefore its subspace
	$A=[0,u_A]$ is Tychonoff.
\end{proof}

We next observe that the topological lifting yields
a stronger structural property: the unit interval
is canonically a retract of the ambient topological
$\ell$-group.

\begin{proposition}\label{prop:canonical-retract}
	Every topological MV-algebra is a retract of a topological
	$\ell$-group. More precisely, under the identification
	$A=[0,u_A]\subseteq G_A$, the map
	$r_A:G_A\to A$ defined by $r_A(g)=g^+\wedge u_A$ is a continuous
	retraction.
\end{proposition}

\begin{proof}
	Since $G_A$ is a topological $\ell$-group, $r_A$ is continuous, and
	$r_A(a)=a$ for every $a\in[0,u_A]$.
\end{proof}

Proposition \ref{prop:canonical-retract} allows results on retral spaces to be
applied directly to topological MV-algebras.

\begin{corollary}\label{cor:sigma-pseudocompact-Knaster}
	Every $\sigma$-pseudocompact topological MV-algebra has
	the Knaster property, and hence is ccc.
\end{corollary}

\begin{proof}
	Let $(A,\tau)$ be a $\sigma$-pseudocompact topological MV-algebra.
	By Proposition~\ref{prop:canonical-retract}, $A$ is a retral space.
	Since $A$ is also a $\Sigma^{-}(\aleph_0)$-space,
	Theorem~\ref{thm:feng-gartside-knaster} implies that $A$ has
	the Knaster property, and hence is ccc by
	Proposition~\ref{prop:knaster-implies-ccc}.
\end{proof}

A \emph{Mal'tsev space} is a topological space admitting
a continuous ternary operation $p:X^3\to X$ satisfying
$p(x,y,y)=x=p(y,y,x)$.
Thus Mal'tsev spaces can be viewed as topological
algebras with a ternary operation.
In particular, topological MV-algebras form an important
class of Mal'tsev spaces
\cite[Proposition~3.1]{XieYang2026}.

Feng and Gartside~\cite[Question~4.16(i)]{FengGartside2026}
asked whether every $\sigma$-pseudocompact Mal'tsev space
is ccc and whether it is retral.
Proposition~\ref{prop:canonical-retract} and
Corollary~\ref{cor:sigma-pseudocompact-Knaster}
provide affirmative answers to both questions
within the class of topological MV-algebras.

Moreover, Corollary~\ref{cor:sigma-pseudocompact-Knaster}
strengthens the ccc result of Xie and
Yang~\cite[Theorem~3.6]{XieYang2026}
by weakening $\sigma$-compactness to
$\sigma$-pseudocompactness and strengthening
ccc to the Knaster property.

We conclude this subsection by applying the canonical
lifting to study paracompactness under local
$\sigma$-compactness.

\begin{proposition}\label{prop:locally-sigma-compact-paracompact}
	Every locally $\sigma$-compact $T_0$ topological MV-algebra
	is paracompact and normal. In particular, every locally
	compact $T_0$ topological MV-algebra is paracompact and normal.
\end{proposition}

\begin{proof}
	Let $(A,\tau)$ be locally $\sigma$-compact and $T_0$.
	Then $G_A$ is Hausdorff by
	Proposition~\ref{prop:canonical-transfer}.
	Choose a $\sigma$-compact neighbourhood $U$ of $0$ in $A$
	and put $B=U-U$. Since subtraction is continuous,
	$B$ is $\sigma$-compact. Moreover, $U^\triangle\subseteq B$
	by the decomposition $g=g^+-g^-$.
	
	By the definition of the canonical topology,
	$W=\sum_{i=1}^{\infty}U^\triangle$ is a $0$-neighbourhood
	in $G_A$, and $W\subseteq\bigcup_{n\geq1}nB$,
	where $nB$ denotes the $n$-fold sum of $B$.
	Since addition is continuous, each $nB$ is $\sigma$-compact,
	and hence so is $\bigcup_{n\ge1}nB$.
	As this set contains the $0$-neighbourhood $W$,
	the group $G_A$ is locally $\sigma$-compact.
	By Proposition~\ref{prop:locally-sigma-compact-group},
	$G_A$ is paracompact.
	
	Finally, $A=[0,u_A]=G_A^+\cap(u_A-G_A^+)$ is closed
	in $G_A$, and hence is paracompact.
	Since $A$ is Hausdorff, it is also normal by
	Proposition~\ref{prop:paracompact-hausdorff-normal}.
\end{proof}

\subsection{Riesz MV-algebras and Koopman representation}
\label{subsec:riesz-koopman}

We next extend the canonical topological lifting
to Riesz MV-algebras
(Definition~\ref{def:riesz-mv-algebra})
and their associated Riesz spaces
(Definition~\ref{def:riesz-space}).
We then identify the lifted topology with a mixed
topology under suitable assumptions.
Finally, by lifting strongly continuous semigroups,
we obtain a unit-interval Koopman representation
for continuous dynamical systems.

A \emph{topological Riesz MV-algebra} is a Riesz MV-algebra
$(A,\tau)$ whose MV-reduct is a topological MV-algebra
and whose scalar multiplication
$[0,1]\times A\to A$ is jointly continuous.
A \emph{topological Riesz space} is a Riesz space
equipped with a topology making vector addition,
scalar multiplication, and lattice operations continuous.
No Hausdorff or local solidity assumption is imposed.

Let $\mathcal T\!\mathit{op}\mathcal{RMV}$ denote the category
of topological Riesz MV-algebras and continuous Riesz
MV-homomorphisms, and let
$\mathcal T\!\mathit{op}\,u\mathcal{R}\mathit{iesz}$
denote the category of topological Riesz spaces with
strong unit and continuous unit-preserving Riesz homomorphisms.

\begin{proposition}
	\label{prop:topological-riesz-lifting}
	Let $(A,\tau)$ be a topological Riesz MV-algebra.
	Then its canonical lifting
	$\Xi_R(A,\tau)=(E_A,u_A,\tau_{E_A})$
	is a topological Riesz space.
	
	Consequently, the canonical functors induce an adjunction
	\[
	\Xi_R\dashv\Gamma_R:
	\mathcal T\!\mathit{op}\mathcal{RMV}
	\rightharpoonup
	\mathcal T\!\mathit{op}\,u\mathcal{R}\mathit{iesz},
	\]
	where $\Xi_R$ is fully faithful.
	
	Moreover, $\tau_{E_A}$ is locally solid if and
	only if $A$ has a zero-neighbourhood basis
	consisting of lower sets.
\end{proposition}

\begin{proof}
	By Theorem~\ref{thm:topological-lifting},
	$(E_A,\tau_{E_A})$ is a topological $\ell$-group
	whose topology induces $\tau$ on $A=[0,u_A]$.
	It remains to prove that scalar multiplication is
	jointly continuous.
	
	We first establish equicontinuity at zero
	for scalar multiplication by elements of $[-1,1]$.
	Let
	\[
	W=\sum_{n=1}^{\infty}U_n^\triangle
	\]
	be a basic zero-neighbourhood of $E_A$.
	For each $n$ and $r\in[0,1]$, joint continuity
	of scalar multiplication on $A$ at $(r,0)$
	yields neighbourhoods $I_r$ of $r$ and
	$V_{n,r}$ of $0$ such that
	$I_rV_{n,r}\subseteq U_n$.
	By compactness of $[0,1]$, finitely many
	$I_{r_1},I_{r_2},\ldots,I_{r_k}$ cover $[0,1]$.
	Setting $V_n=\bigcap_{j=1}^k V_{n,r_j}$,
	we obtain $[0,1]V_n\subseteq U_n$.
	
	If $|\lambda|\leq1$ and $g\in V_n^\triangle$,
	then $(\lambda g)^+,(\lambda g)^-\in U_n$.
	Indeed, for $\lambda\geq0$,
	$(\lambda g)^\pm=\lambda g^\pm$, while for
	$\lambda<0$ the positive and negative parts
	are interchanged. Consequently,
	\[
	[-1,1]V\subseteq W,
	\qquad
	V=\sum_{n=1}^{\infty}V_n^\triangle.
	\]
	
	This extends to every bounded scalar interval.
	Let $M>0$ and let $W$ be a zero-neighbourhood
	of $E_A$. 
	Choose an integer $m\geq M$
	and a zero-neighbourhood $W_0$ such that
	\[
	\underbrace{W_0+W_0+\cdots+W_0}_{m\text{ times}}
	\subseteq W.
	\]
	By the preceding argument, there exists a
	zero-neighbourhood $V$ satisfying
	$[-1,1]V\subseteq W_0$.
	For every $|\lambda|\leq M$, write
	$\lambda=m\mu$, where $|\mu|\leq1$.
	Then, for every $x\in V$,
	\[
	\lambda x
	=\underbrace{\mu x+\mu x+\cdots+\mu x}_{m\text{ times}}
	\in \underbrace{W_0+W_0+\cdots+W_0}_{m\text{ times}}
	\subseteq W.
	\]
	Thus the family of maps
	$\{x\mapsto\lambda x:|\lambda|\leq M\}$
	is equicontinuous at zero, in the sense that
	the same zero-neighbourhood $V$ can be chosen
	for all $|\lambda|\leq M$.
	
	We next establish continuity of the map
	$\lambda\mapsto\lambda x$ for each fixed $x\in E_A$.
	By the strong-unit property and
	Lemma~\ref{lem:layer-decomposition}, we may write
	\[
	x=a_1+a_2+\cdots+a_p-b_1-b_2-\cdots-b_q,
	\qquad a_i,b_j\in A.
	\]
	For each $a_i,b_j\in A$, the maps
	$\lambda\mapsto\lambda a_i$ and
	$\lambda\mapsto\lambda b_j$ from $[0,1]$
	to $E_A$ are continuous.
	Since addition and inversion are continuous,
	the map
	$f_x:[0,1]\to E_A$, defined by
	$f_x(\lambda)=\lambda x$,
	is continuous at zero.
	
	Now let $W$ be any zero-neighbourhood in $E_A$.
	Choose a symmetric zero-neighbourhood $W_0\subseteq W$.
	There exists $\delta>0$ such that
	$f_x(\lambda)\in W_0$ whenever
	$0\leq\lambda<\delta$.
	Since $f_x(-\lambda)=-f_x(\lambda)$,
	we obtain $\lambda x\in W$ whenever
	$|\lambda|<\delta$.
	Thus $\lambda\mapsto\lambda x$ is continuous at zero.
	Since it is a homomorphism of additive topological
	groups, it is continuous on all of $\mathbb R$.
	
	Finally, we prove joint continuity. Fix
	$(\lambda_0,x_0)\in\mathbb R\times E_A$
	and let $W$ be a zero-neighbourhood in $E_A$.
	Choose zero-neighbourhoods $W_1,W_2$ such that
	$W_1+W_2\subseteq W$.
	
	Take $M>|\lambda_0|+1$. By the preceding
	equicontinuity property, there exists a
	zero-neighbourhood $V$ such that $[-M,M]V\subseteq W_1.$
	Since the map $\lambda\mapsto\lambda x_0$
	is continuous at zero, there exists
	$\delta\in(0,1)$ such that
	$|t|<\delta
	\quad\Longrightarrow\quad
	tx_0\in W_2.$
	Now, whenever $|\lambda-\lambda_0|<\delta$
	and $x-x_0\in V$, we have
	\[
	\begin{aligned}
		\lambda x-\lambda_0x_0
		&=\lambda(x-x_0)+(\lambda-\lambda_0)x_0\\
		&\in W_1+W_2\subseteq W.
	\end{aligned}
	\]
	
	Thus scalar multiplication
	$\mathbb R\times E_A\to E_A$
	is jointly continuous, making
	$(E_A,\tau_{E_A})$ a topological Riesz space.
	
	We now verify that the canonical constructions define
	functors between the corresponding topological Riesz
	categories. The functor $\Gamma_R$ is well defined,
	since the operations and homomorphisms restrict
	continuously to the unit intervals.
	
	For $\Xi_R$, the assertion on objects has already
	been established. Let $f:A\to B$ be a continuous
	Riesz MV-homomorphism. By
	Theorem~\ref{thm:riesz-mv-equivalence},
	its canonical extension $\Xi_R(f):E_A\to E_B$
	is a unit-preserving Riesz homomorphism.
	Since its underlying $\ell$-group homomorphism
	is the canonical extension of the MV-reduct of $f$,
	Theorem~\ref{thm:categorical-topological-lifting}
	implies that $\Xi_R(f)$ is continuous.
	Functoriality follows from the canonical lifting.
	
	The unit and counit of the topological MV-adjunction
	are continuous, and by
	Theorem~\ref{thm:riesz-mv-equivalence},
	they are respectively a Riesz MV-isomorphism
	and a Riesz homomorphism.
	Hence they restrict to the topological Riesz
	categories, where the triangle identities
	remain valid.
	Hence
	$\Xi_R\dashv\Gamma_R$.
	Since its unit is a natural isomorphism,
	$\Xi_R$ is fully faithful.
	
	It remains to establish the local-solidity criterion.
	Suppose that $A$ has a zero-neighbourhood
	basis consisting of lower sets.
	Let
	\[
	W=\sum_{n=1}^{\infty}U_n^\triangle
	\]
	be a basic zero-neighbourhood of $E_A$,
	where each $U_n$ is lower.
	
	By continuity of $\oplus$, choose zero-neighbourhoods
	$V_n$ such that $V_n\oplus V_n\subseteq U_n$.
	Set
	\[
	S_n=\{g\in E_A:|g|\in U_n\}.
	\]
	Since $U_n$ is lower, $S_n$ is solid.
	For $g\in V_n^\triangle$, we have
	$|g|=g^+\vee g^-\leq u_A$, and hence
	\[
	|g|=g^++g^-=(g^++g^-)\wedge u_A
	=g^+\oplus g^-\in U_n.
	\]
	Thus $V_n^\triangle\subseteq S_n$.
	Moreover, $S_n\subseteq U_n^\triangle$
	follows from $g^\pm\leq|g|$
	and the lowerness of $U_n$.
	
	We first show that the sum of two solid sets is solid.
	Let $S,T\subseteq E_A$ be solid, and suppose that
	$s\in S$, $t\in T$, and $|z|\leq|s+t|$.
	Put $u=|s|$ and define
	\[
	a=(z\vee(-u))\wedge u,
	\qquad b=z-a.
	\]
	Clearly, $|a|\leq u$.
	Using the lattice identities
	$x\vee y=x+(y-x)^+$ and
	$x\wedge y=x-(x-y)^+$, we obtain
	\[
	\begin{aligned}
		a
		&=(z\vee(-u))
		-\bigl((z\vee(-u))-u\bigr)^+\\
		&=z+(-z-u)^+
		-\bigl((z-u)\vee(-2u)\bigr)^+\\
		&=z+(-z-u)^+-(z-u)^+,
	\end{aligned}
	\]
	where the last equality follows from $u\geq0$.
	Consequently,
	$b=(z-u)^+-(-z-u)^+$.
	Since
	\[
	\begin{aligned}
		(z-u)^+\wedge(-z-u)^+
		&=\bigl((z-u)\wedge(-z-u)\bigr)^+\\
		&=(-|z|-u)^+=0,
	\end{aligned}
	\]
	the two positive terms are disjoint. Hence
	\[
	\begin{aligned}
		|b|
		&=(z-u)^++(-z-u)^+\\
		&=(z-u)^+\vee(-z-u)^+\\
		&=\bigl((z-u)\vee(-z-u)\bigr)^+\\
		&=\bigl((z\vee(-z))-u\bigr)^+\\
		&=(|z|-u)^+\leq|t|,
	\end{aligned}
	\]
	where the last inequality follows from
	$|z|\leq|s+t|\leq u+|t|$.
	Thus $a\in S$ and $b\in T$, so
	$z=a+b\in S+T$.
	By induction, every finite sum of solid sets is solid.
	Since each $S_n$ contains
	$0$, the finite sums $S_1+S_2+\cdots+S_m$ form an
	increasing sequence of solid sets. Hence
	\[
	S=\sum_{n=1}^{\infty}S_n
	\]
	is solid. Moreover, the inclusions
	$V_n^\triangle\subseteq S_n\subseteq U_n^\triangle$
	yield
	\[
	\sum_{n=1}^{\infty}V_n^\triangle
	\subseteq S\subseteq W.
	\]
	Since the left-hand side is a zero-neighbourhood
	in the canonical topology, $S$ is a solid
	zero-neighbourhood contained in $W$.
	Thus $\tau_{E_A}$ is locally solid.
	
	Conversely, if $\tau_{E_A}$ is locally solid,
	intersecting its solid zero-neighbourhood basis
	with $A=[0,u_A]$ yields a zero-neighbourhood
	basis of lower sets in $A$.
\end{proof}

We next identify the canonical topology with
Wiweger's mixed topology under suitable assumptions.

\begin{proposition}
	\label{prop:canonical-mixedification}
	Let $E$ be a Riesz space with strong unit $u$, and let $\sigma$ be
	a Hausdorff topology making $E$ a topological Riesz space.
	Suppose that $\sigma\subseteq\nu_u$,
	where $\nu_u$ is the topology induced by
	the order-unit norm \eqref{eq:order-unit-norm}.
	Then, under the canonical algebraic identification,
	\[
	\Xi_R\Gamma_R(E,u,\sigma)
	=
	\bigl(E,u,\gamma(\nu_u,\sigma)\bigr).
	\]
	
	In particular, the conclusion holds whenever $\sigma$ is
	Hausdorff and locally solid.
\end{proposition}

\begin{proof}
	Every Hausdorff topological Riesz space is
	Archimedean. Indeed, its positive cone is closed,
	since $x\mapsto x^-$ is continuous.
	If $nx\leq y$ for all $n\geq1$, then
	$\frac{1}{n}y-x\geq0$.
	As $\frac{1}{n}y-x\to -x$, the closedness
	of the positive cone yields $-x\geq0$.
	Thus the order-unit norm $\|\cdot\|_u$
	is well defined.
	
	Put $A=[0,u]$ and $\tau=\sigma|_A$.
	Let $\tau_E$ denote the canonical lifting
	of $\tau$.
	By Remark~\ref{rem:canonical-topology}(1),
	$\tau_E$ is the finest group topology
	inducing $\tau$ on $A$.
	
	We claim that $\sigma$ and $\tau_E$
	agree on every order-bounded subset of $E$.
	Indeed, by Lemma~\ref{lem:layer-decomposition},
	for each $n\geq1$, there are maps
	\[
	L_n:[-nu,nu]\longrightarrow A^{2n},
	\qquad
	R_n:A^{2n}\longrightarrow E,
	\]
	such that
	$R_n\circ L_n=\operatorname{id}_{[-nu,nu]}$.
	These maps are continuous with respect
	to both $\sigma$ and $\tau_E$.
	
	Since $\sigma|_A=\tau_E|_A$,
	the corresponding product topologies on
	$A^{2n}$ coincide.
	Moreover, for each of the two topologies,
	$L_n$ is a topological embedding of
	$[-nu,nu]$ into $A^{2n}$, whose inverse
	on its image is the restriction of $R_n$.
	Thus both $\sigma|_{[-nu,nu]}$ and
	$\tau_E|_{[-nu,nu]}$ are transported by
	the same map $L_n$ to the same subspace
	topology on $L_n([-nu,nu])$.
	Hence
	\[
	\sigma|_{[-nu,nu]}=\tau_E|_{[-nu,nu]}.
	\]
	Since $u$ is a strong unit, every order-bounded
	subset of $E$ is contained in $[-nu,nu]$
	for some $n$. Thus $\sigma$ and $\tau_E$
	agree on every order-bounded subset,
	and in particular on every $\nu_u$-bounded subset.

	By Remark~\ref{rem:canonical-topology}(1),
	we have $\sigma\subseteq\tau_E$.
	Moreover, since $\sigma\subseteq\nu_u$
	and the two topologies agree on $[-u,u]$,
	every $\tau_E$-zero-neighbourhood contains
	a sufficiently small order-unit norm ball,
	which lies in $[-u,u]$.
	Hence
	\[
	\sigma\subseteq\tau_E\subseteq\nu_u.
	\]
	Since $\tau_E$ agrees with $\sigma$ on every
	$\nu_u$-bounded subset, the defining property
	of Wiweger's mixed topology yields
	$\tau_E\subseteq\gamma(\nu_u,\sigma)$.
	
	Conversely, since $A$ is $\nu_u$-bounded,
	\[
	\gamma(\nu_u,\sigma)|_A
	=\sigma|_A=\tau.
	\]
	By Remark~\ref{rem:canonical-topology}(1),
	we obtain
	$\gamma(\nu_u,\sigma)\subseteq\tau_E$.
	Therefore
	\[
	\tau_E=\gamma(\nu_u,\sigma).
	\]

	Finally, if $\sigma$ is locally solid,
	every solid zero-neighbourhood contains
	$\varepsilon u$ for some $\varepsilon>0$,
	and hence contains the order-unit norm ball
	of radius $\varepsilon$.
	Thus $\sigma\subseteq\nu_u$.
\end{proof}

Thus, the canonical lifting of the unit-interval topology
recovers Wiweger's mixed topology under the above assumptions.

We now apply the preceding results to the theory
of continuous-time dynamical systems.
A \emph{strongly continuous one-parameter semigroup}
of continuous endomorphisms of a topological algebra $A$
is a family $(S(t))_{t\geq0}$ satisfying
$S(0)=\operatorname{id}_A$,
$S(t+s)=S(t)\circ S(s)$ for all $t,s\geq0$,
and the continuity of $t\mapsto S(t)x$
for every $x\in A$
(cf. \cite[Definition~I.5.1]{EngelNagel2000}).

\begin{proposition}[Unit-interval semigroup correspondence]
	\label{prop:mv-semigroup-lifting}
	Let $(E,u,\sigma)$ be a Hausdorff topological Riesz
	space with strong unit. Put $A=[0,u]$ and
	$\tau=\sigma|_A$. Suppose that $\sigma$ coincides
	with the canonical topology on $E$ obtained by
	lifting $\tau$ from $A$.
	Then there is a one-to-one correspondence between
	\begin{enumerate}
		\item semigroups $(T(t))_{t\geq0}$ of continuous
		unital Riesz endomorphisms of $(E,u,\sigma)$;
		\item semigroups $(S(t))_{t\geq0}$ of continuous
		MV-endomorphisms of $(A,\tau)$.
	\end{enumerate}
	No compatibility with scalar multiplication is
	assumed for the MV-endomorphisms in {\rm(2)}.
	
	The correspondence is given by
	$S(t)=T(t)|_A$, while in the reverse direction
	$T(t)$ is the unique unital Riesz endomorphism
	of $E$ extending $S(t)$.
	
	Moreover, under this correspondence,
	$(T(t))_{t\geq0}$ is strongly continuous if and
	only if $(S(t))_{t\geq0}$ is strongly continuous.
	
	In particular, if $\rho$ is a Hausdorff topology
	making $E$ a topological Riesz space and satisfying
	$\rho\subseteq\nu_u$, then the above correspondence
	holds between semigroups of continuous unital Riesz
	endomorphisms of
	$\bigl(E,u,\gamma(\nu_u,\rho)\bigr)$
	and semigroups of continuous MV-endomorphisms of
	$\bigl([0,u],\rho|_{[0,u]}\bigr)$.
\end{proposition}

\begin{proof}
	We first show that every continuous MV-endomorphism
	$f:A\to A$ preserves scalar multiplication.
	By Theorem~\ref{thm:categorical-topological-lifting}
	and the hypothesis on $\sigma$, $f$ extends
	uniquely to a continuous unital $\ell$-group
	endomorphism $\widetilde f:E\to E$.
	
	Since the additive group of $E$ is uniquely
	divisible, $\widetilde f$ preserves rational
	scalar multiplication. For any $x\in E$ and
	$r\in\mathbb R$, choose a sequence of rational
	numbers $q_n\to r$. By continuity,
	\[
	\begin{aligned}
		\widetilde f(rx)
		&=\lim_{n\to\infty}\widetilde f(q_nx)\\
		&=\lim_{n\to\infty}q_n\widetilde f(x)
		=r\widetilde f(x).
	\end{aligned}
	\]
	Here the Hausdorff assumption ensures uniqueness
	of limits. Thus $\widetilde f$ is real-linear
	and hence a continuous unital Riesz endomorphism.
	In particular, $f$ preserves scalar multiplication
	on $A$.
	
	Now let $(S(t))_{t\geq0}$ be a semigroup of
	continuous MV-endomorphisms of $A$.
	By the preceding argument, each $S(t)$ extends
	uniquely to a continuous unital Riesz endomorphism
	$T(t)$ of $E$.
	Since $S(t+s)=S(t)\circ S(s)$,
	both $T(t+s)$ and $T(t)\circ T(s)$
	are unital Riesz endomorphisms extending
	$S(t+s)$. Similarly, $T(0)$ and
	$\operatorname{id}_E$ both extend
	$S(0)=\operatorname{id}_A$.
	By uniqueness of the extensions,
	\[
	T(t+s)=T(t)\circ T(s),
	\qquad T(0)=\operatorname{id}_E.
	\]
	
	Thus $(T(t))_{t\geq0}$ is a semigroup of
	continuous unital Riesz endomorphisms of $E$.
	
	Conversely, let $(T(t))_{t\geq0}$ be a semigroup
	of continuous unital Riesz endomorphisms of $E$.
	Since each $T(t)$ is positive and preserves $u$,
	we have $T(t)(A)\subseteq A$.
	By Theorem~\ref{thm:riesz-mv-equivalence},
	each restriction $S(t)=T(t)|_A$ is a
	Riesz MV-endomorphism, and hence an MV-endomorphism.
	It is continuous since $\tau=\sigma|_A$,
	and the semigroup identities pass to the restrictions.
	The restriction and extension constructions
	are inverse to each other by uniqueness
	of the extensions, establishing the claimed
	one-to-one correspondence.
	
	It remains to compare strong continuity.
	Suppose first that $(S(t))_{t\geq0}$ is
	strongly continuous. Fix $x\in E$ and choose
	$n\geq1$ such that $|x|\leq nu$.
	By Lemma~\ref{lem:layer-decomposition},
	\[
	x=\sum_{k=1}^nP_k(x)-\sum_{k=1}^nN_k(x),
	\qquad P_k(x),N_k(x)\in A.
	\]
	Since $T(t)|_A=S(t)$, we obtain
	\[
	T(t)x
	=\sum_{k=1}^nS(t)P_k(x)
	-\sum_{k=1}^nS(t)N_k(x).
	\]
	For each $k$, the maps
	$t\mapsto S(t)P_k(x)$ and
	$t\mapsto S(t)N_k(x)$ are continuous
	as $A$-valued maps, and hence as $E$-valued
	maps because $\tau=\sigma|_A$.
	By continuity of addition and subtraction
	in $(E,\sigma)$, the map $t\mapsto T(t)x$
	is continuous. Thus $(T(t))_{t\geq0}$
	is strongly continuous.
	
	Conversely, if $(T(t))_{t\geq0}$ is strongly
	continuous, then for every $a\in A$,
	the map $t\mapsto T(t)a=S(t)a$ is continuous
	as an $E$-valued map taking values in $A$.
	Since $\tau=\sigma|_A$, it is also continuous
	as an $A$-valued map. Hence
	$(S(t))_{t\geq0}$ is strongly continuous.
	
	The final mixed-topology assertion follows from
	Proposition~\ref{prop:canonical-mixedification}.
\end{proof}

We now specialize to spaces of bounded continuous
functions and their associated Koopman semigroups.
The Koopman approach, originating with
Koopman~\cite{Koopman1931}, studies dynamical
systems through their induced linear operators
on spaces of functions. For a jointly continuous
semiflow $\varphi=(\varphi_t)_{t\geq0}$ on a
topological space $X$, the associated Koopman
semigroup acts by
\[
\mathcal K_\varphi(t)f=f\circ\varphi_t,
\qquad t\geq0.
\]

Let $X$ be a nonempty completely regular
Hausdorff space. Consider $E_X=C_b(X,\mathbb R)$
under the pointwise order. It is a Riesz space
with strong unit $u=\mathbf1$, where $\mathbf1$
denotes the constant function with value $1$.
The corresponding Riesz MV-algebra is
\[
A_X=[0,\mathbf1]=C_b(X,[0,1]).
\]

Let $\tau_c^{\mathbb R}$ denote the compact-open
topology on $E_X$, that is, the locally convex
topology generated by the family of seminorms
$\{p_K:K\subseteq X\text{ compact}\}$
(see \cite{FarkasKreidler2020}), where
\[
p_K(f)=\sup_{x\in K}|f(x)|.
\]
Explicitly, a zero-neighbourhood basis consists
of the sets
\[
\bigcap_{j=1}^m
\{f\in E_X:p_{K_j}(f)<\varepsilon_j\},
\]
where $m\geq1$, each $K_j\subseteq X$ is compact,
and each $\varepsilon_j>0$.
The compact-open topology $\tau_c$ on $A_X$
is precisely $\tau_c^{\mathbb R}|_{A_X}$.
Similarly, $\tau_c^{\mathbb C}$ denotes
the compact-open topology on $C_b(X,\mathbb C)$.

The order-unit norm associated with $\mathbf1$
is the supremum norm
\[
\|f\|_u=\|f\|_\infty
=\sup_{x\in X}|f(x)|.
\]
Since $p_K(f)\leq\|f\|_\infty$ for every
compact $K\subseteq X$, we have
$\tau_c^{\mathbb R}\subseteq\nu_u$.

The topology $\tau_c^{\mathbb R}$ is Hausdorff,
since the seminorms $p_K$ separate points.
Moreover, the inequality
$p_K(|f|-|g|)\leq p_K(f-g)$
shows that the modulus map is continuous.
Consequently, the lattice operations are
continuous, and $(E_X,\tau_c^{\mathbb R})$
is a Hausdorff topological Riesz space.

Let $\nu_\infty^{\mathbb C}$ denote the
supremum-norm topology on $C_b(X,\mathbb C)$.
We write
\begin{equation}\label{eq:cb-mixed-topologies}
	\begin{aligned}
		\tau_m^{\mathbb R}
		&=\gamma(\nu_u,\tau_c^{\mathbb R}),\\
		\tau_m^{\mathbb C}
		&=\gamma(\nu_\infty^{\mathbb C},\tau_c^{\mathbb C})
	\end{aligned}
\end{equation}
for the corresponding mixed topologies.
By Proposition~\ref{prop:canonical-mixedification},
the canonical lifting of $(A_X,\tau_c)$ is
\[
\Xi_R(A_X,\tau_c)
=(E_X,\mathbf1,\tau_m^{\mathbb R}).
\]
The complex mixed topology $\tau_m^{\mathbb C}$
is precisely the topology considered by
Farkas and Kreidler
\cite[Definition~2.6]{FarkasKreidler2020}.

Recall that $X$ is \emph{compactly generated},
or a \emph{$k$-space}, if a subset of $X$
is closed whenever its intersection with every
compact subset is relatively closed.
Farkas and Kreidler proved that, when $X$
is compactly generated, a strongly continuous
semigroup $(T(t))_{t\geq0}$ of continuous
complex-linear operators on
$\bigl(C_b(X,\mathbb C),\tau_m^{\mathbb C}\bigr)$
is induced by a unique jointly continuous
semiflow on $X$ if and only if
\[
T(t)\mathbf1=\mathbf1,
\qquad
T(t)|f|=|T(t)f|
\]
for every $f\in C_b(X,\mathbb C)$ and
$t\geq0$; see
\cite[Theorem~3.1]{FarkasKreidler2020}.

This result extends the classical Koopman
representation for compact Hausdorff spaces.
Indeed, when $X$ is compact, we have
$C_b(X,\mathbb C)=C(X,\mathbb C)$ and
$\tau_c^{\mathbb C}=\nu_\infty^{\mathbb C}$,
so the mixed topology coincides with the
supremum-norm topology.
Thus the classical representation is recovered
as a special case, while the mixed topology
allows the same characterization to extend
to compactly generated completely regular spaces.

Combining this representation with
Propositions~\ref{prop:canonical-mixedification}
and~\ref{prop:mv-semigroup-lifting},
we obtain the following unit-interval version.

\begin{theorem}[Unit-interval Koopman representation]
	\label{thm:unit-interval-koopman}
	Let $X$ be a nonempty compactly generated
	completely regular Hausdorff space.
	Then the following classes are in
	one-to-one correspondence:
	\begin{enumerate}
		\item jointly continuous semiflows
		$(\varphi_t)_{t\geq0}$ on $X$;
		
		\item strongly continuous semigroups
		$(T(t))_{t\geq0}$ of continuous
		complex-linear operators on
		$\bigl(C_b(X,\mathbb C),
		\tau_m^{\mathbb C}\bigr)$ such that
		\[
		T(t)\mathbf1=\mathbf1,
		\qquad
		T(t)|f|=|T(t)f|
		\]
		for all $f\in C_b(X,\mathbb C)$
		and $t\geq0$;
		
		\item strongly continuous semigroups
		$(S(t))_{t\geq0}$ of continuous
		MV-endomorphisms of $(A_X,\tau_c)$.
	\end{enumerate}
	Under these correspondences, the semigroup
	associated with a semiflow
	$(\varphi_t)_{t\geq0}$ in {\rm(3)}
	is given by
	\[
	S(t)f=f\circ\varphi_t,
	\qquad f\in A_X.
	\]
\end{theorem}

\begin{proof}
	The correspondence between {\rm(1)}
	and {\rm(2)} is the Koopman representation
	theorem of Farkas and Kreidler
	\cite[Theorem~3.1]{FarkasKreidler2020}.
	It remains to establish the correspondence
	between {\rm(2)} and {\rm(3)}.
	
	We first compare the real and complex
	mixed topologies defined in
	\eqref{eq:cb-mixed-topologies}.
	By the weighted seminorm descriptions
	of mixed topologies
	\cite[p.~475]{Nowak2016}, both
	$\tau_m^{\mathbb R}$ and
	$\tau_m^{\mathbb C}$ are generated by
	seminorms of the form
	\[
	p_w(h)=\sup_{x\in X}w(x)|h(x)|,
	\]
	where $w:X\to[0,\infty)$ is bounded
	and vanishes at infinity, meaning that
	for every $\varepsilon>0$ there exists
	a compact $K\subseteq X$ such that
	$w(x)<\varepsilon$ outside $K$.
	For $f,g\in E_X$, we have
	\[
	\max\{p_w(f),p_w(g)\}
	\leq p_w(f+ig)
	\leq2\max\{p_w(f),p_w(g)\}.
	\]
	Consequently, the identification
	\[
	C_b(X,\mathbb C)\cong E_X^2,
	\qquad
	f+ig\longleftrightarrow(f,g),
	\]
	is a homeomorphism for
	$\tau_m^{\mathbb C}$ and
	$\tau_m^{\mathbb R}\times\tau_m^{\mathbb R}$.
	In particular,
	\[
	\tau_m^{\mathbb C}|_{E_X}
	=\tau_m^{\mathbb R}.
	\]
	
	Suppose first that $(T(t))_{t\geq0}$
	is as in {\rm(2)}.
	For every nonnegative $f\in E_X$,
	we have
	\[
	|T(t)f|=T(t)|f|=T(t)f,
	\]
	so $T(t)f$ is real-valued and nonnegative.
	Since every real-valued function is a
	difference of two nonnegative functions,
	$E_X$ is invariant under $T(t)$.
	Thus $R(t):=T(t)|_{E_X}$ is a
	real-linear, unital, and modulus-preserving
	operator, and hence a unital Riesz
	endomorphism of $E_X$.
	
	By the real-complex mixed-topology
	identification established above,
	each $R(t)$ is continuous for
	$\tau_m^{\mathbb R}$.
	Moreover, for every $f\in E_X$,
	the orbit $t\mapsto R(t)f=T(t)f$
	is continuous in $\tau_m^{\mathbb R}$.
	The semigroup identities pass to the
	restrictions, so $(R(t))_{t\geq0}$
	is a strongly continuous semigroup
	of continuous unital Riesz endomorphisms
	of $(E_X,\tau_m^{\mathbb R})$.
	
	By Proposition~\ref{prop:canonical-mixedification},
	$\tau_m^{\mathbb R}$ is the canonical
	topology lifted from $(A_X,\tau_c)$.
	Hence Proposition~\ref{prop:mv-semigroup-lifting}
	implies that $(S(t))_{t\geq0}$, where
	$S(t)=R(t)|_{A_X}$, is a strongly
	continuous semigroup of continuous
	MV-endomorphisms of $(A_X,\tau_c)$.
	This gives the correspondence
	from {\rm(2)} to {\rm(3)}.
	
	Conversely, let $(S(t))_{t\geq0}$
	be as in {\rm(3)}.
	By Propositions~\ref{prop:canonical-mixedification}
	and~\ref{prop:mv-semigroup-lifting},
	it extends uniquely to a strongly
	continuous semigroup $(R(t))_{t\geq0}$
	of continuous unital Riesz endomorphisms
	of $(E_X,\tau_m^{\mathbb R})$.
	
	For each $t\geq0$, define the
	complexification of $R(t)$ by
	\[
	T(t)(f+ig)=R(t)f+iR(t)g,
	\qquad f,g\in E_X.
	\]
	Since $R(t)$ is real-linear,
	$T(t)$ is complex-linear.
	We now verify the modulus-preservation
	condition in {\rm(2)}.
	
	Fix $t\geq0$ and write $R=R(t)$.
	Since $R$ is positive and unital,
	\[
	|Rh|=R|h|
	\leq\|h\|_\infty R\mathbf1
	=\|h\|_\infty\mathbf1
	\]
	for every $h\in E_X$.
	Thus $R$ is a contraction for the
	supremum norm.
	
	For each $n\geq1$, let
	\[
	F_n=\{2k\pi/n:k=0,1,\ldots,n-1\}
	\]
	and define
	\[
	p_n(a,b)
	=\max_{\theta\in F_n}
	(a\cos\theta+b\sin\theta).
	\]
	As $n\to\infty$, the functions
	$p_n(a,b)$ converge to
	$\sqrt{a^2+b^2}$ uniformly on
	bounded subsets of $\mathbb R^2$.
	Since $R$ is real-linear and
	preserves finite joins,
	\[
	R\bigl(p_n(f,g)\bigr)
	=p_n(Rf,Rg)
	\]
	for all $f,g\in E_X$.
	Passing to the supremum-norm limit,
	using the contractivity of $R$, yields
	\[
	R\sqrt{f^2+g^2}
	=\sqrt{(Rf)^2+(Rg)^2}.
	\]
	Consequently, for $z=f+ig$,
	\[
	T(t)|z|
	=R(t)|z|
	=|R(t)f+iR(t)g|
	=|T(t)z|.
	\]
	Moreover, $T(t)\mathbf1=\mathbf1$.
	Thus $T(t)$ satisfies the
	modulus-preservation and unitality
	conditions in {\rm(2)}.
	
	It remains to verify continuity and
	the semigroup identities.
	Under the real-complex mixed-topology
	identification established above,
	$T(t)$ corresponds to
	$R(t)\times R(t)$.
	Since $R(t)$ is continuous for
	$\tau_m^{\mathbb R}$, each $T(t)$
	is continuous for $\tau_m^{\mathbb C}$.
	Furthermore, for fixed $f,g\in E_X$,
	the orbit
	$t\mapsto T(t)(f+ig)=R(t)f+iR(t)g$
	is continuous, since both
	$t\mapsto R(t)f$ and
	$t\mapsto R(t)g$ are continuous
	in $\tau_m^{\mathbb R}$.
	The semigroup identities for
	$(T(t))_{t\geq0}$ follow directly
	from those of $(R(t))_{t\geq0}$.
	Therefore $(T(t))_{t\geq0}$
	belongs to {\rm(2)}.
	
	Finally, the two constructions are inverse.
	Starting from {\rm(3)}, extending each
	$S(t)$ to $R(t)$ and then complexifying
	recovers $S(t)$ upon restriction to $A_X$.
	Conversely, starting from {\rm(2)},
	the restriction $R(t)=T(t)|_{E_X}$
	is the unique unital Riesz extension
	of $S(t)=R(t)|_{A_X}$ by
	Proposition~\ref{prop:mv-semigroup-lifting}.
	Since $T(t)$ is complex-linear,
	complexifying $R(t)$ recovers $T(t)$.
	Hence {\rm(2)} and {\rm(3)} are
	in one-to-one correspondence,
	completing the proof.
\end{proof}




\section{Conclusion and future directions}
\label{sec:conclusion}

Building on our earlier work~\cite{JiaYang2026}, we have resolved the topological lifting problem for Mundici equivalence in full generality. We prove that every topological MV-algebra admits a canonical representation as the unit interval of a topological $\ell$-group (Theorem~\ref{thm:topological-lifting}), and establish the corresponding categorical adjunction (Theorem~\ref{thm:categorical-topological-lifting}).

The resulting framework yields new topological consequences for MV-algebras, including a group-theoretic proof of the equivalence of $T_0$ and $T_{3\frac12}$ (Corollary~\ref{cor:separation-equivalence}) and the canonical retral-space representation (Proposition~\ref{prop:canonical-retract}). It also extends to topological Riesz MV-algebras (Proposition~\ref{prop:topological-riesz-lifting}), establishes a connection with mixed and strict topologies (Proposition~\ref{prop:canonical-mixedification}), and provides a unit-interval formulation of Koopman representations (Theorem~\ref{thm:unit-interval-koopman}).

Several directions deserve further investigation. First, it is natural to ask whether the topological lifting framework can be extended to pseudo-MV-algebras~\cite{Dvurecenskij2002,GeorgescuIorgulescu2001}, generalized MV-algebras~\cite{GalatosTsinakis2005}, and other related algebraic structures. Second, a further direction is to exploit the lifting framework to derive new results both within topological MV-algebra theory and in related areas of mathematics. Finally, the unit-interval formulation of Koopman theory suggests a systematic study of dynamical systems through the algebraic and topological properties of their associated MV-endomorphism semigroups.

\end{document}